\documentclass[review,onefignum,onetabnum,final]{siamart171218}

\usepackage{lipsum}
\usepackage{amsfonts}
\usepackage{graphicx}
\usepackage{epstopdf}
\usepackage{algorithmic}
\ifpdf
  \DeclareGraphicsExtensions{.eps,.pdf,.png,.jpg}
\else
  \DeclareGraphicsExtensions{.eps}
\fi

\newsiamremark{remark}{Remark}
\newsiamremark{hypothesis}{Hypothesis}
\crefname{hypothesis}{Hypothesis}{Hypotheses}
\newsiamthm{claim}{Claim}

\headers{Spatial profiles of a diffusive epidemic model}{P. Rui, R. Salako, and Y. Wu}

\title{Spatial profiles of  a reaction-diffusion epidemic model with nonlinear incidence mechanism and constant total population \thanks{Submitted to the editors DATE.
\funding{R. Peng was partially supported by NSF of China (No. 12271486, 12171176).}}}

\author{Rui Peng\thanks{School of Mathematical Sciences, Zhejiang Normal
University, Jinhua, 321004, Zhejiang, China 
  (\email{pengrui\_seu@163.com}).}
\and Rachidi Salako\thanks{Department of Mathematical Sciences, University of Nevada Las Vegas, Las Vegas, NV 89154, USA 
  (\email{rachidi.salako@unlv.edu}).}
\and Yixiang Wu\thanks{ Department of Mathematical Sciences, Middle Tennessee State University, Murfreesboro, TN 37132, USA 
  (\email{yixiang.wu@mtsu.edu}).}}
\usepackage{amsopn}

\usepackage{soul}
\usepackage{graphics}
\usepackage{caption}
\usepackage{subcaption}

\newtheorem{tm}{Theorem}[section]
\newtheorem{prop}{Proposition}[section]

\newtheorem{lem}{Lemma}[section]

\newtheorem{rk}{Remark}[section]

\numberwithin{equation}{section}
\numberwithin{tm}{section}

\def\p{\partial}
\def\R{{\mathbb R}}

\def\be{\beta}

\def\De{\Delta}
\def\de{\delta}
\def\ga{\gamma}

\def\Om{\Omega}

\usepackage{placeins}

\begin{document}

\maketitle

\begin{abstract}
  This paper examines a susceptible-infected-susceptible  (SIS) epidemic reaction-diffusion model with no-flux boundary conditions and constant total population. The infection mechanism in the model is described by a nonlinear term of the form $S^qI^p$ with $0<p\leq1$ and $q>0$.
  We explore the spatial profiles of the endemic equilibrium with respect to  small movement rates of susceptible and/or infected populations. Our results extend and improve existing results in the literature, providing further insights into how transmission mechanisms and population mobility limitations impact the spread of infectious diseases. We reinforce and complement the theoretical results with numerical simulations.
\end{abstract}



\section{Introduction}

\subsection{Background}
The outbreaks of new and existing infectious diseases are facilitated by various factors, such as human or animal invasions into new ecosystems, global warming, environmental degradation, increased international travel, and changes in economic patterns \cite{AM,AM2,DH,Martcheva}. Mathematical models and computer simulations have emerged as valuable tools to analyze the spread and control of infectious diseases. They enable the construction and testing of theories, evaluation of quantitative conjectures, addressing specific questions, assessing sensitivities to parameter changes, and estimation of key parameters from data. Understanding the transmission dynamics of infectious diseases at the community, regional, and national levels can lead to more effective strategies to reduce the transmission or eradicate these diseases \cite{BCC,Het2,McCallum}.

It is commonly known that transmission mechanisms play an essential role in modeling the infectious
diseases. In the pioneering work of Kermack and McKendrick \cite{KM1}, the disease transmission
mechanism was chosen to be the so-called mass action type: the number of newly infected individuals
per unit area per unit of time is given by a bilinear incidence function $\beta SI$, where $S$ and $I$ are the density of susceptible and infected individuals, respectively, and $\beta$ is the disease transmission rate. As pointed out in several works (e.g.,  \cite{DV1,Het1,McCallum}),
this mechanism may not be effective enough to describe the evolution process of certain diseases. Moreover, to better describe the interactions of susceptible and infected people, the transmission incidence function should take into consideration of the population characteristics and behavior  that capture the essential features of the disease. One widely used incidence function is the power-law incidence $\beta S^qI^p$, where $p$ and $q$ are positive constants (e.g., see \cite{Het2,HLV,HV,LHL,LLI,LR} and the references therein).

On the other hand, spatial heterogeneity and mobility are important factors that can greatly affect the transmission and spread of infectious diseases, which are usually neglected in classical ordinary differential equation (ODE) for epidemic models.   Here, spatial heterogeneity refers to the variation in physical, social, and environmental conditions across space, which can lead to differentiations in disease transmissions. For example, urban areas may have higher population densities and more social interactions compared to rural areas resulting in larger disease transmission rates. Spatial mobility refers to the movement of individuals across space, which can facilitate the spread of communicable diseases. For example, individuals who travel from one region to another can introduce the disease into new areas and contribute to the spatial expansion of the disease; the movement of infected individuals within a region can create clusters of infections and affect the spatial distribution of the disease.

In accordance with this perspective, reaction-diffusion epidemic models with varying coefficients have been proposed and studied to comprehend the combined influence of spatial heterogeneity and mobility on disease transmissions (\cite{allen2008asymptotic,LPW,Ruan}).
In this two-paper series, we focus  on a  specific class of  spatial epidemic models, namely the reaction-diffusion SIS epidemic models with nonlinear incidence function $\beta S^qI^p$. The chosen incidence function offers the flexibility to model disease transmission comprehensively, encompassing a wide spectrum of dynamics and including the previously mentioned mass action  mechanism (i.e., $p=q=1$). Recently,  there are many related works on diffusive models with mass action or standard incidence ($\beta SI/(S+I)$) mechanisms  \cite{allen2008asymptotic,ChenCui,Cui,CLPZ,CL,CLL,DP,GLPZ2024,KMP,LLL,LXZ,LPW,LPX,LT2023,LLou,LS,LS2023,LSS,MWW1,MWW2,Peng,PL,Sa,SLX,TM,WWK,WW}. Our work contributes to this growing field by investigating  models with nonlinear incidence mechanism $\beta S^qI^p$.

\subsection{The model} The  SIS epidemic model under consideration in this paper reads as follows:
\begin{equation}\label{model-1}
    \begin{cases}
        \p_tS=d_S\Delta S-\beta(x) S^qI^p+\gamma(x) I,\ \ \ & x\in\Om,  t>0,\cr
        \p_t I= d_I\Delta I +\beta(x) S^qI^p-\gamma(x) I, & x\in \Om,\ t>0,\cr
        \p_{\nu}S=\p_{\nu}I=0, & x\in\p \Om, \ t>0,\cr
        S(x,0)=S_0(x),\ \ I(x,0)=I_0(x),&x\in\Omega,
    \end{cases}
\end{equation}
where $\Omega$ is a bounded domain in $\R^n$ with  smooth boundary $\p\Om$; $\nu$ is the outward unit normal to $\p\Om$; and $\beta$ and $\gamma$ are positive H\"older continuous functions on $\bar{\Om}$, representing disease transmission and recovery rates, respectively. The two unknown functions $S(x,t)$ and $I(x,t)$ are, respectively, the density of  susceptible and infected individuals at  spatial location $x$ and time $t$. The homogeneous Neumann boundary condition means that there is no population flux across the boundary $\p\Om$. The positive constants $d_S$ and $d_I$ are the movement rate of susceptible and infected individuals, respectively. Integrating the sum of the first two equations of \eqref{model-1} over $\Omega$, one has
\begin{equation}\label{model-aa}
\int_\Omega(S(x,t)+I(x,t))\,{\rm d}x=\int_\Omega(S_0(x)+I_0(x))\,{\rm d}x=:N,\ \ \ \forall t\ge 0,
\end{equation}
which means that the total population $N$ is constant. Throughout this paper, we always assume that $N$ is a fixed positive constant.
One may refer to \cite{FCGT,Het2,PW2021} for more biological background of \eqref{model-1}.

An equilibrium of  \eqref{model-1}-\eqref{model-aa} is a classical solution of the elliptic system:
\begin{equation}\label{1-1}
\begin{cases}
d_S\Delta S-\beta(x) S^qI^p+\gamma(x) I =0, \ \ \ &  x\in\Om, \cr
d_I\Delta I +\beta(x) S^qI^p-\gamma(x) I=0, & x\in \Om,\cr
\p_{\nu}S=\p_{\nu}I=0, & x\in\p \Om, \cr
\displaystyle
\int_{\Om}(I+S)=N.
\end{cases}
\end{equation}
A nonnegative equilibrium $(S, I)$ is called a \textit{disease free equilibrium} (DFE) if $I=0$ and an \textit{endemic equilibrium} (EE) if $I\neq 0$. Clearly,
$E_0:=(N/|\Omega|,0)$ is the unique DFE; any EE $(S,I)$ of \eqref{1-1}  satisfies $S,\,I>0$ on $\bar\Omega$ due to the well-known maximum principle for elliptic equations. In the current paper, we only focus on  the parameter range of $0<p\leq 1$ and $q>0$. Models with $p > 1$ can exhibit bistable dynamics \cite{FCGT,PW2021}, characterized by the coexistence of one stable positive equilibrium and one stable disease-free equilibrium. Such Allee effect type phenomenon will not be discussed here.

The main objective of this paper is to investigate the asymptotic behavior of the EE of  system  \eqref{1-1} as the diffusion rates $d_S$ and/or $d_I$ tend to zero. These profiles can provide valuable insights into understanding the geographic distribution of the disease and predicting the potential outcomes of disease control strategies that limit population movement.
For system  \eqref{1-1}, previous studies have focused on the case $q=p=1$ \cite{castellano2022effect,deng2016dynamics,deng2023dynamics,PWZZ2023,wu2016asymptotic,WJL2018} and the case $0<p<1$ in one spatial domain \cite{li2020asymptotic}. Our findings presented in this paper significantly improve and complement the previous research. In doing so, we have to develop new techniques to overcome the substantial challenges associated with working in higher dimensional domains or unexplored scenarios in the existing literature.

The current work is aimed at examining the spatial profiles of solutions in SIS-type reaction-diffusion epidemic models characterized by the nonlinear incidence mechanisms mentioned above. In a companion paper \cite{PSW2024b}, we investigate the spatial profiles of solutions to the following SIS diffusive epidemic model:
\begin{equation}\label{model-2}
\begin{cases}
\p_tS=d_S\De S+\Lambda(x)-S-\beta(x) S^qI^p+\gamma(x) I,\ \ \ & x\in\Omega, t>0,\cr
\p_tI=d_I\De I +\beta(x) S^qI^p-(\gamma(x)+\eta(x)) I, & x\in\Omega, t>0,\cr
\p_{\nu}S=\partial_{\nu}I=0, & x\in\p\Omega, t>0, \cr
 S(x,0)=S_0(x),\ \ I(x,0)=I_0(x),&x\in\Omega.
\end{cases}
\end{equation}
In \eqref{model-2}, $\Lambda-S$ is  the recruitment term, representing that the susceptible population is subject to a linear growth (the coefficient before $S$ is assumed to be a positive constant and normalized to one for simplicity of presentation); $\eta$ stands for disease induced mortality rate; $\beta$ and $\gamma$  have the same meaning as in \eqref{model-1}. All coefficient functions  are positive H\"older continuous functions on $\bar{\Om}$. Note that the total population in system \eqref{model-2} varies as opposed to  system \eqref{model-1}. For detailed biological interpretations of \eqref{model-2}, one may see \cite{Het1,Het2,li2018diffusive}.

It is worth mentioning that most of the techniques employed in the analysis of \eqref{model-1} do not apply to system \eqref{model-2}.  In fact, further new ideas and techniques are developed in \cite{PSW2024b} to explore the spatial profiles of solutions of system \eqref{model-2}. A comprehensive comparison of the findings for the two epidemic models can be found in \cite{PSW2024b}.

The rest of this paper is  structured as follows. Section \ref{Main-Results-section} presents the main results, and their proofs are provided in Sections \ref{case p=1} and \ref{case p<1}. Section 3 contains some preliminary results. In Section \ref{Nuremical-Sim-section}, we discuss the practical implications of our findings in terms of disease control measures, and we also present numerical simulations focusing on the case where $p=1$.


\section{Main results}\label{Main-Results-section}
In this section, we will present the main results regarding the asymptotic profiles of the EE as $d_S$ and/or $d_I$ approach zero for systems \eqref{1-1}.

The value of $p$ has a notable influence on the dynamics of model \eqref{model-1}. Specifically, the asymptotic profiles of the EEs in relation to the population diffusion rates differ when $0<p<1$ compared to the case of $p=1$. For $0<p<1$, Theorem \ref{theorem_ex}-{\rm (i)} provides results on the well-posedness and eventual boundedness of classical solutions for system \eqref{model-1}+\eqref{model-aa}.
When $p=1$,  as in \cite{allen2008asymptotic,PW2021,PZ2012}), one can define the basic reproduction number $\mathcal{R}_0$  as
\begin{equation}\label{R_0-q}
\mathcal{R}_{0}=   \sup_{\varphi\in W^{1,2}(\Om)\setminus\{0\}}\frac{\Big(\frac{N}{|\Omega|}\Big)^q\int_{\Om}\be\varphi^2}{\int_{\Om}(d_I|\nabla \varphi|^2+\gamma\varphi^2)},
 \end{equation}
which is independent of $d_S$.   It is well known that the supremum in the definition of  $\mathcal{R}_0$ is achieved \cite{allen2008asymptotic}. We recall some  basic properties of $\mathcal{R}_0$ in Lemma \ref{lem_R0}. As in most epidemic models, Theorem \ref{theorem_ex}-{\rm (ii)} confirms that when $p=1$, $\mathcal{R}_0>1$ is enough for the persistence of the disease and existence of EE solutions of system \eqref{model-1}+\eqref{model-aa}.

For convenience,  define
\begin{equation}\label{r-def}
r(x)=\frac{\ga(x)}{\be(x)},\ \ \quad \forall\, x\in\bar{\Om}.
\end{equation}
We call $r$ the risk function for  \eqref{1-1}, which characterizes the risk of each point on $\bar\Omega$ when $d_S=d_I=0$ ($1/r$ is  proportional to the basic reproduction number of the corresponding ODE model  when $p=1$).

Given a continuous function $f$ on $\bar\Omega$, we denote
 $$
 f_{\text{min}}=\min_{x\in\bar\Omega}f(x)\ \ \ \text{and} \ \ f_{\text{max}}=\max_{x\in\bar\Omega}f(x).
 $$

From now on, unless otherwise stated, we will denote  $(S, I)$  as an  EE  of  \eqref{1-1}  if it exists.

\subsection{Asymptotic profile of EE as $d_I\to 0$} Define
\begin{equation}\nonumber
\Om_{*}:=\Big\{x\in\bar{\Om}:\ r(x)=r_{\min} \Big\},
\end{equation}
which is the set of highest-risk locations. We have the following result regarding the limit of $(S, I)$ as $d_I\to 0$.

\begin{tm}\label{TH2.1} Suppose that $d_S,\ q>0$. Then the following statements hold.
\begin{enumerate}
    \item[\rm (i)] Suppose that $p=1$. Let  $\frac{N}{|\Om|}>r_{\min}^{\frac{1}{q}}$ such that  \eqref{1-1} has at least one EE for $0<d_I\ll 1$. Then we have
\begin{equation}\label{Eq8-01}
    \lim_{d_I\to0}\|S-r_{\min}^{\frac{1}{q}}\|_{L^{\infty}(\Omega)}=0, \quad \lim_{d_I\to 0}\int_{\Om}I=N-|\Omega|r_{\min}^{\frac{1}{q}},
\end{equation}
and
$$
I\to0 \quad \text{locally uniformly in} \  \bar\Omega\setminus\Omega_* \ as \  d_I\to0.
$$
Moreover, if $\Om_*$ consists of finitely many isolated points, say $\Om_*=\{x_j\in\bar\Omega:\ 1\leq j\leq m\}$, then  there exist nonnegative numbers $c_1,\cdots, c_m$ satisfying 
$$
\sum_{j=1}^mc_j=N-|\Omega|r_{\min}^{\frac{1}{q}},
$$
such that up to a subsequence,
$$
I\to \sum_{j=1}^mc_j\delta_{x_j} \  \text{weakly-star in}\  [C(\bar\Omega)]^*\ as \ d_I\to 0,
$$
where $\delta_{x_j}$ is the Dirac measure centered at $x_j$.

\item[\rm (ii)] Suppose that $0<p<1$.  Then \eqref{1-1} has at least one EE  for any $d_I>0$, and
\begin{equation}\label{S-limit-d_I=0-equ2}
(S,I)\to\left(S_*,\ \left(\frac{S_*^{q}}{r}\right)^{\frac{1}{1-p}}\right)\ \ \mbox{uniformly on}\ \bar\Omega \ \ \mbox{as}\ d_I\to0,
\end{equation}
where the constant $S_*$ is the unique positive solution of the algebraic equation
\begin{equation}\label{kappa-star-equ}
    N=\int_{\Om}\Big[S_* +\Big(\frac{S_*^{q}}{r}\Big)^{\frac{1}{1-p}}\Big].
\end{equation}

\end{enumerate}
\end{tm}

 For a subset $A\subset \R^n$, we denote by $\mathring{A}$ its interior, that is the largest open subset of $A$. If $p=1$ and $|\Omega_*|\neq 0$, Theorem \ref{TH2.1}(i) can be improved as follows.
\begin{prop} \label{Prop2.1}
Suppose that $d_S>0$, $p=q=1$, and  $\frac{N}{|\Om|}>r_{\min}$. If $\mathring{\Omega}_*\ne\emptyset$, then there exists $d_0>0$ independent of $d_S$ such that \eqref{1-1} has a unique EE $(S, I)$ for any $0<d_I\le \min\{d_S,d_0\}$. Moreover, there exists $C>0$ depending on $d_S$ such that
\begin{equation}\label{TH2.2-eq1}
    \|S-r_{\min}\|_{L^{\infty}(\Omega)}\le C\frac{d_I}{d_S}, \ \quad \Big|\int_{\Om}I-(N-|\Omega|r_{\min})\Big|\le Cd_I|\Omega|,  \ \quad   \|I\|_{\infty}\le C,
\end{equation}
and
$$
I\to0 \  \text{locally uniformly in}\ \ \bar\Omega\setminus\Omega_*\  \text{as}\ d_I\to0.
$$
In addition, restricted to a subsequence if necessary, we have $I\to \hat{I}$  weakly in $W^{1,2}(\Omega)$, where $\hat{I}\in L^{\infty}(\Om)\cap W^{1,2}(\Om)$ satisfies that $\hat I=0$ on $\bar\Omega\setminus{\Om}_*$.

In particular, if  $\Om_*=\cup_{j=0}^m{\Om}_*^j \ (0\leq m<\infty)$
with $|{\Om}_*^0|=0$ and each connected smooth component ${\Om}_*^j \ (1\leq m<\infty)$ fulfilling ${\Om}_*^j\cap{\Om}_*^i=\emptyset$ for all $0\leq i,\,j\leq m$ and $i\not=j$, then
in each ${\Om}_*^j\ (1\leq j\leq m)$, either $\hat{I}=0$ or $\hat{I}\in C^{2}(\mathring{\Om}_*^j)$ is a positive.  Furthermore, setting $ \Upsilon:=\{1\leq i\leq m:\ \hat I>0\ \mbox{in}\ \mathring{\Om}_*^i\}$, there exist $\hat{a}_j>0$, $j\in \Upsilon$, satisfying 
\begin{equation}\label{2-8ab}
 \displaystyle
      \sum_{j\in \Upsilon}\int_{{\Om}_*^j}\hat{I}=N-|\Omega|r_{\min},
\end{equation}
such that for each $j\in \Upsilon$,  $\hat{I}$ is a classical solution of
\begin{equation}\label{2-8aa}
  \begin{cases}
      \Delta \hat{I}+\frac{\beta}{d_S}(\hat{a}_j-\hat{I})\hat{I}=0,\ \ \ & x\in \mathring{\Om}_*^j,\cr
      \partial_{\nu}\hat I=0, &x\in\partial\Om\cap\partial\Om_*^j,\cr
      \hat I=0, &x\in\partial\Om_*^j\setminus\partial\Om.
      \end{cases}
\end{equation}
\end{prop}

\vskip6pt
\begin{rk}\label{r-th2.1} Theorem \ref{TH2.1} and Proposition \ref{Prop2.1} improve/generalize some results in the literature.
\begin{enumerate}
    \item[\rm (1)] If $p=q=1$, Theorem \ref{TH2.1} {\rm(i)} was obtained in \cite{castellano2022effect} and Proposition \ref{Prop2.1} was proved when the domain is one-dimensional in  \cite{PWZZ2023}.

    \item[\rm (2)] If $0<p<1$, under an extra condition of either $q\geq 1$ or $p\geq 1-q>0$, Theorem \ref{TH2.1} {\rm (ii)} was established in \cite{li2020asymptotic}  when the domain $\Omega$ is one-dimensional.
   \end{enumerate}
\end{rk}

One might wonder if the profile of the EE in Theorem \ref{TH2.1}(i) can be obtained by taking $p\to 1^-$ in Theorem \ref{TH2.1}(ii). We provide the following remark about this question.

\begin{rk}
    For     $0<p<1$, let $S_*$ be the unique positive solution of \eqref{kappa-star-equ}. As $p\to 1^{-}$,  $S_*\to \min\{\frac{N}{|\Omega|},r_{\min}^{\frac{1}{q}}\}$   and $(S_*^{q}/r)^{\frac{1}{1-p}}\to 0$  locally uniformly in $\bar\Om\setminus\Omega_*$. Furthermore, the latter limit holds uniformly on $\bar\Om$ if $N/|\Om|<r_{\min}^{\frac{1}{q}}$. Hence, as $p\to 1^{-}$, if $N/|\Om|<r_{\min}^{\frac{1}{q}}$ then $(S_*, (S_*^{q}/r)^{\frac{1}{1-p}})$ approximates the unique DFE $E_0$ of \eqref{1-1};  if $N/|\Omega|>r_{\min}^{\frac{1}{q}}$,   we recover the limiting profiles described by \eqref{Eq8-01}. However, the profiles of the EEs  as described in \eqref{2-8aa} could not be predicted by letting $p\to 1^{-}$ in Theorem \ref{TH2.1} {\rm (ii)}.
\end{rk}

\subsection{Asymptotic profile of EE as $d_S\to 0$}
We have the following result on the asymptotic profile of the EE for the case $d_S\to 0$.
\begin{tm}\label{TH2.3} Fix $d_I,\ q>0$. Then the following statements hold.

\begin{enumerate}  \item[\rm (i)] Suppose that $p=1$. Let $\mathcal{R}_{0}>1$ such that  \eqref{1-1} has at least one EE $(S, I)$.
\begin{enumerate}
    \item[\rm (i-a)] If $N<\int_{\Om}r^{\frac{1}{q}}$, then there is $C>0$ such that
\begin{equation}\label{Eq8-03-1}
    \|I\|_{L^\infty(\Omega)}\le Cd_S,\ \ \quad \forall\, 0<d_S<1.
\end{equation}
Moreover, restricted to a subsequence if necessary,
 \begin{equation}\label{Eq8-03-02}
    S\to S^*:=\frac{N(1-d_II^*)}{\int_{\Om}(1-d_II^*)}\ \  \ \text{uniformly on}\ \bar{\Om} \ \text{as}\ d_S\to 0,
\end{equation}
where $0<{I}^*<\frac{1}{d_I}$ is a classical solution of
\begin{equation}\label{Eq8-03-03}
    \begin{cases}
        d_I\Delta I^* +\beta\Big(N^q\frac{(1-d_II^*)^q}{\|1-d_II^*\|_{L^1(\Om)}^q}-r\Big)I^*=0,\ \ \ & x\in\Om, \cr
        \p_{\nu} I^*=0, & x\in\p\Om.
    \end{cases}
\end{equation}

\item[\rm (i-b)] If $N=\int_{\Om}r^{\frac{1}{q}}$, then  $I\to0$ uniformly on $\bar\Omega$ as $d_S\to0$.
Moreover, restricted to a subsequence if necessary, either $S\to r^{\frac{1}{q}}$ or $S\to S^*$ uniformly on $\bar\Omega$ as $d_S\to 0$, where $S^*$ is given as in \eqref{Eq8-03-02}.

\item[\rm (i-c)] If $N>\int_{\Om}r^{\frac{1}{q}}$, then there is a sequence of positive numbers $\{d_{S,n}\}_{n\ge1}$ satisfying $d_{S,n}\to0$ as $n\to\infty$, such that for every $n\geq 1$, \eqref{1-1} has an EE solution $(S_n,I_n)$ for $d_S=d_{S,n}$ fulfilling
\begin{equation}\nonumber
(S_n,I_n)\to\left(r^{\frac{1}{q}},\ \frac{1}{|\Om|}\Big(N-\int_{\Om}r^{\frac{1}{q}}\Big)\right)\ \ \mbox{uniformly on}\ \bar\Omega \ \mbox{as}\ n\to\infty.
\end{equation}
\end{enumerate}

\item[\rm (ii)] Suppose that $0<p<1$.  Then \eqref{1-1} has at least one EE $(S, I)$, and it satisfies
\begin{equation}\label{SIc}
 (S,I)\to\left(r^\frac{1}{q} I_*^\frac{1-p}{q},\ I_*\right)\ \ \mbox{uniformly on}\ \bar\Omega \ \mbox{as}\ d_S\to0,
\end{equation}
where the constant $I_*$ is the unique positive solution of the algebraic equation
\begin{equation}\label{Istr}
    N=\int_{\Om}\Big[ \Big(rI_*^{1-p}\Big)^{\frac{1}{q}}+I_*\Big].
\end{equation}

\end{enumerate}
\end{tm}

There are parameter ranges for which each scenario in Theorem \ref{TH2.3}-{\rm (i)} holds, and we show this in the following remark.   
\begin{rk}
For convenience, define 
$$
\mathcal{R}_1 := \sup_{\varphi \in W^{1,2}(\Omega) \setminus \{0\}} \frac{\int_{\Omega} \beta \varphi^2}{\int_{\Omega} (d_I |\nabla \varphi|^2 + \gamma \varphi^2)}.
$$
Then by \eqref{R_0-q}, we have  $ \mathcal{R}_0 = \left( \frac{N}{|\Omega|} \right)^q \mathcal{R}_1$, which implies that $\mathcal{R}_0 > 1$ if and only if  $N > |\Omega| \mathcal{R}_1^{-\frac{1}{q}}$. Suppose $r$ is not constant. By Lemma \ref{lem_R0}, $|\Om|\mathcal{R}_1^{-\frac{1}{q}}$ is strictly increasing in $d_I > 0$ with
$$
\lim_{d_I \to 0^+} |\Om|\mathcal{R}_1^{-\frac{1}{q}} =
|\Om|{r^{\frac{1}{q}}_{\min}} \quad \text{and} \quad \lim_{d_I \to \infty} |\Om|\mathcal{R}_1^{-\frac{1}{q}} = |\Om|\left(\frac{\int_{\Omega} r\beta}{\int_{\Omega} \beta}\right)^{\frac{1}{q}}.
$$
Moreover, using $N$ as a varying parameter, 
 the following conclusions hold:
\begin{itemize}
    \item[\rm(1)]  If $|\Omega|r_{\min}^{\frac{1}{q}}<N<\min\left\{\int_{\Omega}r^{\frac{1}{q}},|\Omega|\left(\frac{\int_{\Omega} r\beta}{\int_{\Omega} \beta}\right)^{\frac{1}{q}}\right\}$, there is $d_*>0$ such that $\mathcal{R}_0>1$ if and only if $0<d_I<d_*$. Therefore,  the conditions of Theorem \ref{TH2.3}-\textit{(i-a)} are satisfied for every $0<d_I<d_*$. 
    
    \item[\rm(2)] If $N=\int_{\Omega}r^{\frac{1}{q}}<|\Omega|\left(\frac{\int_{\Omega}r\beta}{\int_{\Omega} \beta}\right)^{\frac{1}{q}}$,  there is $d_*>0$ such that $\mathcal{R}_0>1$  if and only if $0<d_I<d_*$. Therefore,  the conditions of Theorem \ref{TH2.3}-\textit{(i-b)} hold for every $0<d_I<d_*$.

\item[\rm (3)] If $\int_{\Omega}r^{\frac{1}{q}}<N<|\Omega|\left(\frac{\int_{\Omega} r\beta}{\int_{\Omega} \beta}\right)^{\frac{1}{q}}$, there is $d_*>0$ such that $\mathcal{R}_0>1$ if and only if $0<d_I<d_*$. Therefore,  the conditions of Theorem \ref{TH2.3}-\textit{(i-c)} are  satisfied for every $0<d_I<d_*$.
    
    \item[\rm (4)] 
    If $\max\left\{|\Omega|\left(\frac{\int_{\Omega} r \beta}{\int_{\Omega} \beta}\right)^{\frac{1}{q}},\int_{\Omega}r^{\frac{1}{q}}\right\}<N$ and $N\neq \int_\Omega r^\frac{1}{q}$, then $\mathcal{R}_0>1$ and  the conditions of Theorem \ref{TH2.3}-\textit{(i-c)} are satisfied for every $d_I>0$ .
    
    \item[\rm (5)] If $|\Omega|\left(\frac{\int_{\Omega} r \beta}{\int_{\Omega} \beta}\right)^{\frac{1}{q}}\le N<\int_{\Omega}r^{\frac{1}{q}}$, then $\mathcal{R}_0>1$ and  the conditions of Theorem \ref{TH2.3}-\textit{(i-a)} are satisfied for every $d_I>0$.

    \item[\rm (6)]  If $|\Omega|\left(\frac{\int_{\Omega} r \beta}{\int_{\Omega} \beta}\right)^{\frac{1}{q}}< N=\int_{\Omega}r^{\frac{1}{q}}$, then $\mathcal{R}_0>1$ and  the conditions of Theorem \ref{TH2.3}-\textit{(i-b)} are satisfied for every $d_I>0$.
\end{itemize}
\end{rk}

\begin{rk} Theorem \ref{TH2.3} improves/generalizes some existing results in the literature.
\begin{enumerate}
 \item[\rm (1)] Theorem \ref{TH2.3} {\rm (i-a)} and {\rm (i-c)} were obtained when $p=q=1$ in \cite{castellano2022effect}, but Theorem \ref{TH2.3} {\rm(i-b)} seems to be new even when $p=q=1$.
 \item[\rm (2)] Theorem \ref{TH2.3} {\rm(ii)} was proved in \cite{li2020asymptotic} when the domain $\Omega$ is one-dimensional.
\end{enumerate}
\end{rk}

\vskip6pt
Some statements in Theorem \ref{TH2.3} can be further improved.
\begin{rk}
\begin{enumerate}
    \item[\rm (1)] For Theorem \ref{TH2.3} {\rm (i-b)}, by Remark \ref{r-th2.3-a}{\rm(1)},
when $\beta$ is constant and $q=1$  there is a sequence of EEs, denoted by $(S_n,I_n)$, such that $(S_n,I_n)\to (r^{\frac{1}{q}},0)$ uniformly on $\bar\Omega$ as $d_{S,n}\to 0$. We suspect that this is also true when $\beta$ is non-constant and $q>0$.

\item[\rm (2)]
 For Theorem \ref{TH2.3} {\rm (i-c)}, if $p=q=1$, system \eqref{1-1} also has EEs with the limiting profiles in Theorem \ref{TH2.3} {\rm (i-a)}   under some additional hypotheses (see \cite{Castellano2023multiplicity}). However,
 suppose in addition that $\frac{N}{|\Omega|}\ge r_{\max}^{\frac{1}{q}}$, by Remark \ref{r-th2.3-a}{\rm (2)}, it holds that
$$
(S,I)\to \left(r^{\frac{1}{q}},\ \frac{1}{|\Om|}\Big(N-\int_{\Om}r^{\frac{1}{q}}\Big)\right)\quad\text{uniformly on}\ \bar\Om \ \text{as}\ d_S\to0.
$$
\end{enumerate}
\end{rk}

The following remark states that the results in Theorem \ref{TH2.3}(i) may not be obtained from Theorem \ref{TH2.3}(ii) by taking $p\to 1^-$.
\begin{rk}
  For $0<p<1$, let $I_*$ be the unique positive solution of the algebraic equation \eqref{Istr}. As $p\to 1^{-}$, if $N\le \int_{\Om}r^{\frac{1}{q}}$ then $I_*\to 0$ and $I_*^{\frac{1-p}{q}}\to N/\int_{\Omega}r^{\frac{1}{q}}$; if $N>\int_{\Omega}r^{\frac{1}{q}}$ then $I_*\to (N-\int_{\Om}r^{\frac{1}{q}})/|\Om|$ and $I_*^{\frac{1-p}{q}}\to 1$. Hence, as $p\to 1^{-}$, the following conclusions hold:

\begin{itemize}
\item[\rm (1)] If $N\le \int_{\Om}r^{\frac{1}{q}}$ and $\mathcal{R}_{0}<1$, the EEs  of \eqref{1-1}  approximate $\big(Nr^{\frac{1}{q}}/\int_{\Om}r^{\frac{1}{q}},0\big)$  after taking $d_S\to 0$ first and $p\to 1^-$, which is not equal to the unique DFE $E_0$ when $p=1$ and $r$ is not constant. This suggests that system \eqref{1-1} may have multiple EEs when $p=1$ and $\mathcal{R}_{0}<1$ for some parameter ranges. This result is known to be true when $p=q=1$; see \cite{Castellano2023multiplicity}.

\item[\rm (2)] If $N<\int_{\Om}r^{\frac{1}{q}}$ and $\mathcal{R}_{0}>1$,  the EEs  of \eqref{1-1}  approximate $\big(Nr^{\frac{1}{q}}/\int_{\Om}r^{\frac{1}{q}},0\big)$  after taking $d_S\to 0$ first and $p\to 1^-$. Note that this limiting profile is different from that given by  Theorem \ref{TH2.3} {\rm (i-a)}.

\item[\rm (3)]  If $N>\int_{\Om}r^{\frac{1}{q}}$ and $\mathcal{R}_{0}>1$,  the profiles of the EEs of \eqref{1-1}  after taking $d_S\to 0$ first and $p\to 1^-$ are identical to the limiting profiles of the EEs when $p=1$ as $d_S\to 0$.
\end{itemize}
\end{rk}

\subsection{Asymptotic profile of EE as $d_S,\,d_I\to 0$}
Finally, we are  concerned with the asymptotic profile of the EE when limiting both $d_S$ and $d_I$.

\begin{tm}\label{TH2.4}
Fix $q>0$. Then the following statements hold.
\begin{enumerate}
    \item[\rm (i)] Suppose that $p=1$ and $\sigma>0$. Let $\frac{N}{|\Omega|}>r^{\frac{1}{q}}_{\min}$ such that  \eqref{1-1} has at least one EE  for any $0<d_I\ll 1$.
     Then we have
\begin{equation}\label{TH2.4-eq1}
(S,I)\to\left(S_{\sigma},\ I_{\sigma}\right)\ \ \mbox{uniformly on}\ \bar\Omega \ \mbox{as}\ d_I+|\frac{d_I}{d_S}-\sigma|\to0,
\end{equation}
where $S_{\sigma}:=\min\{\tilde{\kappa}_{\sigma},r^{\frac{1}{q}}\}$, $I_{\sigma}:=\frac{1}{\sigma}\Big(\tilde{\kappa}_{\sigma}-r^{\frac{1}{q}}\Big)_+$,  and  the constant $\tilde{\kappa}_{\sigma}$ is the unique positive solution of the algebraic equation
\begin{equation}\label{TH2.4-eq2}
   N=\int_{\Omega}\Big[\min\{\tilde{\kappa}_{\sigma},r^{\frac{1}{q}}\}+\frac{1}{\sigma}\Big(\tilde{\kappa}_{\sigma}-r^{\frac{1}{q}}\Big)_+\Big].
\end{equation}
Moreover, the following conclusions about $(\tilde{\kappa}_{\sigma},S_{\sigma},I_{\sigma})$ hold.
\begin{enumerate}
    \item[\rm (i-a)] If $N<\int_{\Om}r^{\frac{1}{q}}$, then $(S_{\sigma},I_{\sigma})\to(\min\{\tilde{\kappa}_{\infty},r^{\frac{1}{q}}\},\ 0)$ uniformly on $\bar\Om$ and $\tilde{\kappa}_{\sigma}\to \tilde{\kappa}_{\infty}$ as $\sigma\to\infty$, where $\tilde{\kappa}_{\infty}$ is the unique positive number in $(r_{\min}^{\frac{1}{q}},r_{\max}^{\frac{1}{q}})$ satisfying $N=\int_{\Omega}\min\{\tilde{\kappa}_{\infty},r^{\frac{1}{q}}\}$.

    \item[\rm (i-b)] If $N=\int_{\Om}r^{\frac{1}{q}}$, then  $(S_{\sigma},I_{\sigma})\to(r^{\frac{1}{q}},\ 0)$ uniformly on $\bar\Om$  and $\tilde{\kappa}_{\sigma}\to r_{\max}^{\frac{1}{q}}$ as $\sigma\to\infty$.

    \item[\rm (i-c)] If $N>\int_{\Om}r^{\frac{1}{q}}$, then there is $\sigma^*\gg 1$ such that $\tilde{\kappa}_{\sigma}=\frac{\sigma}{|\Omega|}(N-\int_{\Om}r^{\frac{1}{q}})+\frac{1}{|\Omega|}\int_{\Omega}r^{\frac{1}{q}}$ for any $\sigma\ge \sigma^*$. Moreover,
    $$
    (S_{\sigma},I_{\sigma})\to\left(r^{\frac{1}{q}},\ \frac{1}{|\Omega|}\Big(N-\int_{\Omega}r^{\frac{1}{q}}\Big)\right)\ \  \text{uniformly on}\ \bar\Om\ \text{as}\ \sigma\to\infty.
    $$

    \item[\rm (i-d)] $\tilde{\kappa}\to r_{\min}^{\frac{1}{q}}$, $S_{\sigma}\to r_{\min}^{\frac{1}{q}}$ uniformly on $\bar\Om$, $I_{\sigma}\to 0$ locally uniformly on $\bar\Om\setminus\Om_*$, and $\int_{\Om}I_{\sigma}\to N-|\Omega|r^{\frac{1}{q}}_{\min}$ as $\sigma\to 0^+$. Moreover if $|\Om_*|>0$, then $I_{\sigma}\to \frac{1}{|\Om_*|}(N-|\Om|r^{\frac{1}{q}}_{\min})$ uniformly on $\bar{\Om}_*$ as $\sigma\to 0 $.
\end{enumerate}

    \item[\rm (ii)] Suppose that $0<p<1$, $q>0$ and  $\sigma>0$. Then \eqref{1-1} has at least one EE, and it satisfies
\begin{equation}\nonumber
(S,I)\to\left(r^{\frac{1}{q}}
    \big[I_{\sigma}(\cdot,\tilde{\kappa}_{\sigma})\big]^{\frac{1-p}{q}},\ I_{\sigma}(\cdot,\tilde{\kappa}_{\sigma})\right)\ \ \mbox{uniformly on}\ \bar\Omega\ \ \mbox{as}\ d_I+|\frac{d_I}{d_S}-\sigma|\to0,
\end{equation}
where $(\Tilde{\kappa}_{\sigma},I_{\sigma}(\cdot,\tilde{\sigma}))$ is uniquely determined by \eqref{LU2-1} and  \eqref{LU2-2}.

\end{enumerate}
\end{tm}

\vskip6pt
\begin{rk}\label{r-th2.4} If $p=q=1$,  \eqref{TH2.4-eq1} and \eqref{TH2.4-eq2}  have been established in \cite{wu2016asymptotic}, and the remaining results in Theorem \ref{TH2.4} are new. 
\end{rk}


\section{Preliminaries}  In this section, we prepare some preliminary results.

\subsection{Two useful lemmas} We start by  recalling a basic convexity property of a $C^2$-function on $\bar\Omega$ at its local extrema as follows.

\begin{lem}[{\cite[Lemma 2.2]{jiang2007priori}, \cite[Lemma 3.1]{peng2013qualitative}}]\label{l3.1} Suppose that $w\in C^2({\bar\Omega})$
and $\partial_{\nu}w=0$ on $\partial\Omega$. The following
statements hold.

{\rm (i)} If $w$ attains a local maximum at $x_1\in\bar\Omega$,
then $\nabla w(x_1)=0$ and $\Delta w(x_1)\leq0$.

{\rm (ii)} If $w$ attains a local minimum at $x_2\in\bar\Omega$,
then $\nabla w(x_2)=0$ and $\Delta w(x_2)\geq0$.
 \end{lem}

\vskip6pt
The following lemma  follows from standard arguments in the literature.

\begin{lem}[{\cite[Lemma 2]{Dancer}, \cite[Proposition 2.2]{DL1997}}]\label{Lem-3} Suppose that $a, b>0$ and $\mathcal{O}\subset \mathbb{R}^n$ is a bounded domain with a smooth boundary. Let $u(\cdot,a,b,\mathcal{O})$ be the unique nonnegative stable classical solution of the Fisher-KPP elliptic system
\begin{equation}\label{Eq3-6}
    \begin{cases}
\displaystyle        \Delta u +\frac{\beta}{b}(a-u)u=0,\ \ \ & x\in\mathcal{O},
        \cr u=0, & x\in\p \mathcal{O}.
    \end{cases}
\end{equation}
Then there exists a constant $a_{\rm low}(b)>0$ such that $u(\cdot,a,b,\mathcal{O})>0$ in $\Omega$ for all $a>a_{\rm low}(b)$ and $u(\cdot,a,b,\mathcal{O})\equiv 0$ for all $0<a\le a_{\rm low}(b)$. Moreover for each $b>0$, the function $(a_{\rm low},\infty)\ni a\mapsto u(\cdot,a,b,\mathcal{O})\in C(\overline{\mathcal{O}})$ is continuous, strictly increasing and satisfies
\begin{equation}\label{Eq3-8}
    \lim_{a\to\infty}\min_{x\in K}u(x,a,b,\mathcal{O})=\infty
\end{equation}
for any nonempty compact set $K\subset\mathcal{O}$.
\end{lem}

\subsection{Global dynamics and existence of EEs of system \eqref{model-1}}\label{R-0-subsection}
In this subsection, we recall some preliminary results about model \eqref{model-1}.
In order to ensure the well-posedness of \eqref{model-1}, we impose the following assumptions on the initial data:
\begin{enumerate}
\item[{\rm(A)}] $S_0$ and $I_0$ are nonnegative continuous functions on $\bar\Omega$. Moreover,
\begin{enumerate}
\item[(i)] The initial value $I_0\geq,\not\equiv0$ on $\bar\Omega$;
\item[(ii)] If $0<q< 1$, $S_0(x)>0$ for all $x\in\bar\Omega$;
\item[(iii)] If $0<p<1$, $I_0(x)>0$ for all $x\in\bar\Omega$.
\end{enumerate}
\end{enumerate}

We recall the following result in \cite{PW2021} concerning the existence, uniqueness and uniform persistence of the solutions, and the existence of positive equilibrium for model \eqref{model-1}.
\begin{tm}[{\cite[Theorem 2.3 and Theorem 2.5]{PW2021}}]\label{theorem_ex} Suppose that $q>0,\ 0< p\le 1$ and {\rm (A)} holds. Then the following statements hold.
\begin{enumerate}
\item[{\rm(i)}] System \eqref{model-1} admits a unique global classical solution $(S, I)$  with
$$
S(x, t), I(x, t)>0\quad \text{for all}\ x\in\bar\Omega \ \text{and}\  t>0,
$$
and there exists $M_\infty>0$ depending on the initial data such that
\begin{equation}\nonumber
\|S(\cdot, t)\|_{L^\infty(\Omega)}, \ \ \|I(\cdot, t)\|_{L^\infty(\Omega)} \le M_\infty, \ \  t\ge 0.
\end{equation}
Moreover, there exists $N_\infty>0$ independent of initial data (but depending  on $N$) such that
\begin{equation}\nonumber
\max\{\limsup_{t\rightarrow\infty}\|S(\cdot, t)\|_{L^\infty(\Omega)},  \ \limsup_{t\rightarrow\infty}\|I(\cdot, t)\|_{L^\infty(\Omega)}\} \le N_\infty.
\end{equation}

\item[{\rm(ii)}]
Suppose in addition that $\mathcal{R}_{0}>1$ if $p=1$.
Then there exists $\epsilon_0>0$ independent of initial data (but depending  on $N$) such that for any solution $(S, I)$ of \eqref{model-1}, we have
\begin{equation}\nonumber
\min\{\liminf_{t\rightarrow\infty} S(x, t), \  \liminf_{t\rightarrow\infty}I(x, t)\} \ge \epsilon_0,
\end{equation}
uniformly for $x\in\bar\Omega$. Moreover, \eqref{model-1} has at least one EE.
\end{enumerate}
\end{tm}

By Theorem \ref{theorem_ex}, if $0<p<1$, model \eqref{model-1} admits at least one EE solution; if $p=1$, model \eqref{model-1} admits at least one EE solution if $\mathcal{R}_0>1$.

The following result about the properties of $\mathcal{R}_0$ can be found in the literature \cite{allen2008asymptotic}.

\begin{lem}\label{lem_R0}
    Let $d_I, q>0$ and $\mathcal{R}_0$ be defined by \eqref{R_0-q}. Then, the following conclusions hold.
  \begin{enumerate}
   \item If $\beta/\gamma$ is constant, then $\mathcal{R}_0=\left(\frac{N}{|\Omega|}\right)^q\frac{\beta}{\gamma}$ for any $d_I>0$;
   \item If  $\beta/\gamma$ is not constant, then $\mathcal{R}_0$ is strictly decreasing in $d_I$ with
   $$
\lim_{d_I\to 0}\mathcal{R}_0=\left(\frac{N}{|\Omega|}\right)^q \left(\frac{\beta}{\gamma}\right)_{\max}\ \ \text{and} \ \  \lim_{d_I\to \infty}\mathcal{R}_0=\left(\frac{N}{|\Omega|}\right)^q \frac{\int_\Omega\beta}{\int_\Omega\gamma}.
   $$
    \end{enumerate}
\end{lem}

The proof of the following result is quite similar to that of \cite[Lemmas 3.1-3.2]{allen2008asymptotic},  so we omit it here.
\begin{lem}\label{Lem-2-1} Fix $d_S, d_I>0$. Then the following statements hold.
\begin{enumerate}
\item[\rm(i)] Let $(S,I)$ be an EE of \eqref{1-1}. Then the function
\begin{equation}\label{Eq6-1}
\kappa= d_SS+d_II,\ \ \quad x\in\Omega
\end{equation}
is a positive constant. Furthermore, defining
\begin{equation}\label{Eq6-2}
    \bar{S}=\frac{S}{\kappa} \quad \text{and}\quad  \bar{I}=\frac{I}{\kappa},
\end{equation}
then $(\kappa,\bar{S},\bar{I})$ satisfies
\begin{equation}\label{Eq6-3}
    \bar{S}=\frac{1}{d_S}(1-d_I\bar{I}),
\end{equation}
\begin{equation}\label{Eq6-4}
    N=\frac{\kappa}{d_S} \int_{\Omega}\left[\big(1-d_I\bar{I}\big)+d_S\bar{I} \right]
\end{equation}
and
\begin{equation}\label{Eq6-5}
    \begin{cases}
\displaystyle    d_I\Delta \bar{I} +\kappa^{p-1}\Big( \beta\left[\frac{\kappa}{d_S}(1-d_I\bar{I})\right]^q-\gamma \kappa^{1-p} \bar{I}^{1-p}\Big)\bar{I}^p=0,\ \ & x\in\Omega,\cr
    \partial_{\nu}\bar{I}=0, & x\in\partial\Omega,\cr
    0<\bar{I}<\frac{1}{d_I}, & x\in\Omega.
    \end{cases}
\end{equation}
\item[\rm(ii)] If $(\kappa,\bar{S},\bar{I})$ solves \eqref{Eq6-3}-\eqref{Eq6-5},  then $(S,I)=(\kappa\bar{S},\kappa\bar{I})$ is an EE of \eqref{1-1}.
\end{enumerate}
\end{lem}

\section{Proofs of the main results: case of $p=1$}\label{case p=1} In this section, we suppose that $p=1$ and prove Theorem \ref{TH2.1}(i), Proposition \ref{Prop2.1}, Theorem \ref{TH2.3}(i) and Theorem \ref{TH2.4}(i).

\subsection{Small $d_I$:\ proof of Theorem \ref{TH2.1}{\rm (i)} and Proposition \ref{Prop2.1}}
 We first treat the case of small dispersal rate $d_I$ and establish Theorem \ref{TH2.1}(i) and Proposition \ref{Prop2.1}.

\begin{proof}[Proof of Theorem \ref{TH2.1}{\rm (i)}] By Lemma \ref{lem_R0}, $\mathcal{R}_{0}$ is decreasing in $d_I$ with
$$
\mathcal{R}_{0}\to \left(\frac{N}{|\Om|}\right)^q\left(\frac{\beta}{\gamma}\right)_{\text{max}}=\left(\frac{N}{|\Om|}\right)^q\frac{1}{r_{\min}}.
$$
Since $\frac{N}{|\Om|}>r_{\min}^{\frac{1}{q}}$, there is $d_0>0$ such that $\mathcal{R}_{0}>1$ for all $d_I\in (0, d_0)$. By Theorem \ref{theorem_ex},  \eqref{1-1} has an EE $(S,I)$ for any $d_I\in (0, d_0)$ and $d_S>0$. For any $0<d_I<\min\{d_0,d_S\}$, let $({\kappa},\bar{S},\bar{I})$ be defined as in Lemma \ref{Lem-2-1}. Observe from \eqref{Eq6-4} that
\begin{equation}\label{Eq8-1}
    \frac{\kappa}{d_S}=\frac{N}{|\Om|+(d_{S}-d_I)\int_{\Om}\bar{I}}<\frac{N}{|\Om|},\ \ \ \quad \forall\;  0<d_I<\min\{d_0,d_S\}.
\end{equation}
By \eqref{Eq6-5},  $\bar{I}$ satisfies
$$
\begin{cases}
\displaystyle    d_I\Delta \bar{I} +\beta\Big[\Big(\frac{\kappa}{d_S} \Big)^q-r  \Big]\bar{I}>0,\ \ \ & x\in\Om,\cr
    \partial_{\nu}\bar{I}=0, & x\in\p\Om.
\end{cases}
$$
Since $\bar{I}>0$ on $\bar{\Om}$, it is easy to see that
\begin{equation}\label{Eq8-2}
    \frac{\kappa}{d_S}>r_{\min}^{\frac{1}{q}}.
\end{equation}

We claim that
    \begin{equation}\label{Eq8-5}
        \lim_{d_I\to 0}\frac{\kappa}{d_S}=r_{\min}^{\frac{1}{q}}.
    \end{equation}
Assume to the contrary that the claim is false. By \eqref{Eq8-1}-\eqref{Eq8-2}, up to a subsequence if necessary, we may assume $\frac{\kappa}{d_S}\to \kappa_*$ as $d_I\to0$ for some $\kappa_*\in(r_{\min}^{\frac{1}{q}},\frac{N}{|\Om|}]$ .

Notice that $\tilde{u}:=d_{I}\bar{I}$ satisfies
    \begin{equation}\label{Eq8-3}
        \begin{cases}
      \displaystyle      d_{I}\Delta\tilde{u}+\beta\frac{\kappa^q}{d_S^q}\Big[ (1-\tilde{u})^q-r{\frac{d_S^q}{\kappa^q}} \Big]\tilde{u}=0,\ \ \ &  x\in\Om, \cr
           \partial_{\nu}\tilde{u}=0, & x\in\p \Om,  \cr
           0<\tilde{u}<1.
        \end{cases}
    \end{equation}
As a result, we can employ the singular perturbation theory (see, for instance, \cite[Lemma 2.4]{du2009effect}) to conclude that
    \begin{equation*}
        \lim_{d_I\to 0}\tilde{u}=\Big(1-\frac{r^{\frac{1}{q}}}{\kappa_*}\Big)_{+}\ \ \ \ \text{uniformly on}\ \bar\Om.
    \end{equation*}
    This implies that
    \begin{equation}\label{Eq8-4}
    \lim_{d_I\to 0}\int_{\Om}\tilde{u}=\int_{\Om}\Big(1-\frac{r^{\frac{1}{q}}}{\kappa_*}\Big)_{+}.
    \end{equation}
    On the other hand, by \eqref{Eq8-2}, we have
    $$
    \int_{\Om}\tilde{u}=d_{I}\int_{\Om}\bar{I}
    =\frac{d_{I}}{\kappa}\int_{\Om}I\le\frac{Nd_{I}}{\kappa}\le \frac{d_{I}N}{d_S r_{\min}^{\frac{1}{q}}} \to 0\ \  \text{as}\ d_I\to 0.
    $$
It then follows from \eqref{Eq8-4} that
$$
\int_{\Om}\Big(1-\frac{r^{\frac{1}{q}}}{\kappa_*}\Big)_{+}=0.
$$
Hence  $\kappa_{*}\le r_{\min}^{\frac{1}{q}}$, which  contradicts the assumption that $\kappa_*>r_{\min}^{\frac{1}{q}}$. This proves \eqref{Eq8-5}.

By $N=\int_{\Om}(S+I)$, to complete the proof of \eqref{Eq8-01}, we only need to show that
\begin{equation}\label{Eq8-6}
    \lim_{d_I\to0}\|S-r_{\min}^{\frac{1}{q}}\|_{L^\infty(\Omega)}=0.
\end{equation}
It is clear from \eqref{Eq6-1} that
$$
S=\frac{\kappa}{d_S}-\frac{d_I}{d_S}I<\frac{\kappa}{d_S},\ \ \quad \forall\; 0<d_I<\min\{d_0,d_S\}.
$$
On the other hand, it follows from the first equation of \eqref{1-1} and Lemma \ref{l3.1}(i) that
$$
S\ge r_{\min}^{\frac{1}{q}},\ \ \quad \forall\; 0<d_I<\min\{d_0,d_S\}.
$$
Hence by \eqref{Eq8-5}, \eqref{Eq8-6} holds.

Given $x\in\R^n$ and $\de>0$, we denote by $\mathcal{B}_{x,\delta}$ the open ball centered at $x$ with radius $\de$. Let $K\subseteq \Om\setminus\Om_*$ be compact. Then there exist $0<\sigma\ll 1$ and $x_1,\cdots,x_n\in K$ such that $K\subset\cup_{i=1}^n\mathcal{B}_{x_i,3\sigma}\subset\cup_{i=1}^n\mathcal{B}_{x_i,6\sigma}\subset\Om\setminus\Om_*$.
Let $i=1,\cdots,n$ be fixed. By \eqref{Eq8-6} and the definition of $\Om_*$, there exist $d_0, \eta_0>0$ such that
\begin{equation}\label{Eq8-21}
\beta S^{q}(x)-\gamma<-\eta_0,\ \ \  \ \forall x\in\mathcal{B}_{x_i,5\sigma}, \; 0<d_I<d_0.
\end{equation}
This implies that $\Delta I>0$ and  $I$ is subharmonic on $\mathcal{B}_{x_i,5\sigma}$ for any $0<d_I<d_0$. Hence, by the maximum principle for elliptic equations (\cite{GT}), there is $\mu_0>0$, independent of $d_I$, such that
\begin{equation}\label{Eq8-22}
\|I\|_{L^{\infty}(\mathcal{B}_{x_i,3\sigma})}\le \mu_0\int_{\mathcal{B}_{x_i,4\sigma}}I,\ \ \quad \forall\;  0<d_I<d_{0}.
\end{equation}
Choose a  smooth function $\varphi\in C^{2}(\bar{\Omega})$ such that $0\le \varphi\le 1$, $\varphi\equiv 1$ on $\mathcal{B}_{x_i,4\sigma}$ and $\varphi\equiv 0$ for $x\notin\mathcal{B}_{x_i,5\sigma}$. Multiplying the second equation of \eqref{1-1} by $\varphi$ and integrating by  parts, we have
\begin{eqnarray*}
    0&=&d_I\int_{\Om}I\Delta \varphi+\int_{\mathcal{B}_{x_i,5\sigma}}(\beta S^q-\gamma)I\varphi\\
    &<&d_I\int_{\Om}I\Delta \varphi-\eta_0\int_{\mathcal{B}_{x_i,5\sigma}}I\varphi \\
    &\le& d_I\int_{\Om}I\Delta \varphi-\eta_0\int_{\mathcal{B}_{x_i,4\sigma}}I, \quad\forall\; 0<d_I<d_0,
\end{eqnarray*}
where we have used \eqref{Eq8-21}. The last inequality together with \eqref{Eq8-22} yields
$$
\|I\|_{L^{\infty}(\mathcal{B}_{x_i,3\sigma})}\le \mu_0\int_{\mathcal{B}_{x_i,4\sigma}}I\le \mu_0\frac{d_I}{\eta_0}\|\Delta\varphi\|_{\infty}\int_{\Om}I\le \mu_0\frac{d_I}{\eta_0}\|\Delta\varphi\|_{\infty}N\to 0,\ \quad \text{as}\ d_{I}\to 0.
$$
Since $i$ is arbitrary and $K\subset\cup_{i=1}^n\mathcal{B}_{x_i,3\sigma}$, we conclude that $I\to 0$ uniformly on $K$ as $d_I\to0$.
For the more general case that $K\subseteq \bar\Om\setminus\Om_*$ is compact, one can easily adapt the above analysis to show  that $I\to 0$ uniformly on $K$ as $d_I\to0$. Since $K\subset \bar{\Omega}\setminus\Omega_*$ is an arbitrary compact set,   $I\to 0$ locally uniformly in $\bar\Om\setminus\Om_*$. 

Finally, suppose that  $\Om_*=\{x_j\in\bar\Omega:\ 1\leq j\leq m\}$. Since $I$ is bounded in $L^1(\Omega)$, $I\to
\sum_{j=1}^m c_j\delta_{x_j}$ weakly in $L^1(\Omega)$ as $d_I\to 0$, where $\delta_{x_j}$ is the Dirac measure centered at $x_j$. Since $\int_{\Om}I\to N-|\Omega|r_{\min}^{\frac{1}{q}}$ as $d_I\to0$, the nonnegative constants $c_j$ satisfy $\sum_{j=1}^mc_j=N-|\Omega|r_{\min}^{\frac{1}{q}}$. This completes the proof of Theorem \ref{TH2.1}(i).
\end{proof}

%
%

We  need the following lemma for the special case $q=p=1$.
\begin{lem}\label{Lem-4} Suppose that $q=p=1$, $\mathring{\Om}_*\ne\emptyset$  and $\frac{N}{|\Om|}>r_{\min}$. Then there exists $d_0>0$ such that $\mathcal{R}_{0}>1$ for all $d_I\in (0, d_0)$. For any $0<d_{I}\le d_0$ and $d_S> 0$, let $(S,I)$ be the EE of \eqref{1-1}. Then for all $d_S>0$, it holds that
\begin{equation}\label{Eq3-9}
    \liminf_{d_I\to0}\frac{d_S}{d_I}\Big(\frac{\kappa}{d_S}-r_{\min}\Big)\ge \frac{N}{|\Omega|}-r_{\min}>0
\end{equation}
and
\begin{equation}\label{Eq3-10}
    \limsup_{d_I\to0}\frac{d_S}{d_I}\Big(\frac{\kappa}{d_S}-r_{\min}\Big)<\infty.
\end{equation}
\end{lem}
\begin{proof}
The existence of $d_0$ and EE can be ensured as in the proof of Theorem \ref{TH2.1}{\rm (i)}. Let $(\kappa,\bar{S},\bar{I})$ be given by Lemma \ref{Lem-2-1}(i).  It follows from \eqref{Eq6-1} that
\begin{equation} \label{Eq3-11}
\frac{d_S}{d_I}\Big(\frac{\kappa}{d_S}-r_{\min}\Big)=\frac{d_S(S-r_{\min})}{d_I}+I.
\end{equation}
Recall that $S\ge r_{\min}$. So we have
$$
\frac{d_S}{d_I}\Big(\frac{\kappa}{d_S}-r_{\min}\Big)\ge \frac{1}{|\Omega|}\int_{\Omega}I.
$$
This combined with \eqref{Eq8-01} proves \eqref{Eq3-9}.

We now proceed to prove \eqref{Eq3-10}. By \eqref{Eq6-5}, $I$ satisfies
   \begin{equation}\label{Eq3-12}
    \begin{cases}
\displaystyle    \Delta I +\frac{\beta}{d_S}\Big[\frac{d_S}{d_I}\Big(\frac{\kappa}{d_S}-r_{\min}\Big)-I\Big]I=0,\ \ \ & x\in\Om_*,\cr
    I>0, & x\in\Om_*.
    \end{cases}
\end{equation}
Since $\mathring{\Om}_* $ is not empty, there exist $x_2\in \Omega$ and $\de_2>0$ such that $\mathcal{B}(x_2,2\delta_2)\subset\Om_*$. Thanks to \eqref{Eq3-12}, one can appeal to a comparison argument and the uniqueness of $u\Big(\cdot,\frac{d_S}{d_I}\Big(\frac{\kappa}{d_S}-r_{\min}\Big),d_S,\mathcal{B}_{x_2,2\delta_2}\Big)$ (see Lemma \ref{Lem-3}) to conclude that
\begin{equation*}
    u\Big(\cdot,\frac{d_S}{d_I}\Big(\frac{\kappa}{d_S}-r_{\min}\Big),d_S,\mathcal{B}_{x_2,2\delta_2}\Big)\le I\ \  \ \mbox{on}\ \bar\Omega,
\end{equation*}
from which we deduce that
\begin{eqnarray*}
&& \min_{x\in\overline{\mathcal{B}}_{x_2,\delta_2} }u\Big(x,\frac{d_S}{d_I}\Big(\frac{\kappa}{d_S}-r_{\min}\Big),d_S,\mathcal{B}_{x_2,2\delta_2}\Big) \\
&\leq&  \frac{1}{|\mathcal{B}_{x_2,\delta_2}|}  \int_{\mathcal{B}_{x_2,\delta_2}}u\Big(\cdot,\frac{d_S}{d_I}\Big(\frac{\kappa}{d_S}-r_{\min}\Big),d_S,\mathcal{B}_{x_2,2\delta_2}\Big) \cr
 &\le& \frac{1}{|\mathcal{B}_{x_2,\delta_2}|}  \int_{\mathcal{B}_{x_2,\delta_2}}I\le \frac{N}{|\mathcal{B}_{x_2,\delta_2}|}.
\end{eqnarray*}
Therefore, in view of Lemma \ref{Lem-3} and \eqref{Eq3-12}, \eqref{Eq3-10} holds.
\end{proof}

We are ready to prove Proposition \ref{Prop2.1}.

\begin{proof}[Proof of Proposition \ref{Prop2.1}] The existence of $d_0$ and EE is guaranteed as in the proof of Theorem \ref{TH2.1}{\rm (i)}. Let $(\kappa,\bar{S},\bar{I})$ be given by Lemma \ref{Lem-2-1}(i).    By \eqref{Eq3-10},
$$
C:=\sup_{0<d_I\le \min\{d_S,d_0\}}\frac{d_S}{d_I}\Big(\frac{\kappa}{d_S}-r_{\min} \Big)<\infty,
$$
which is depending on $d_S$ but independent of $d_I$.
By \eqref{Eq6-1} and $S\ge r_{\min}$, we have
\begin{equation}\label{Eq4-1}
0\le S-r_{\min}< \frac{\kappa}{d_S}-r_{\min}\le C \frac{d_I}{d_S},\ \quad \forall\; d_I\le \min\{d_S,d_0\}.
\end{equation}
This combined with $\int_{\Omega}I=N-\int_{\Omega}S$ proves the first two inequalities of  \eqref{TH2.2-eq1}.
In addition, we observe from \eqref{Eq6-1} that
\begin{equation}\label{Eq4-2}
I=\frac{d_S}{d_I}\Big(\frac{\kappa}{d_S}-r_{\min}\Big)-\frac{d_S}{d_I}(S-r_{\min})\le C,\quad \forall x\in\bar\Omega,\ \, d_I\le\min\{d_S,d_0\}.
\end{equation}
This verifies the last inequality of \eqref{TH2.2-eq1}.

By \eqref{Eq4-1}, we have
\begin{align*}
0 =& \Delta I+\beta\Big( \frac{S-r_{\min}}{d_I}-\frac{r-r_{\min}}{d_I}\Big)I\cr
\le & \Delta I +C\frac{\beta}{d_S}I.
\end{align*}
Multiplying this inequality by $I$ and integrating over $\Omega$, we deduce
$$
\|\nabla I\|_{L^2(\Om)}\le \frac{|\Om| \beta_{\max}  C^2}{d_S},\ \quad \forall\; d_I\le \min\{d_S,d_0\}.
$$
This together with \eqref{Eq4-2} implies that
$$
\sup_{0<d_I\le\min\{d_S,d_0\}}\|I\|_{W^{1,2}(\Om)}<\infty.
$$
Hence passing to a subsequence if necessary, $I\to \hat{I}$ weakly in $W^{1,2}(\Omega)$ and $I\to \hat{I}$ in $L^2(\Omega)$ for some $\hat{I}\in W^{1,2}(\Omega)\cap L^\infty(\Omega)$ as $d_I\to 0$. By Theorem \ref{TH2.1}{\rm (i)}, $I\to 0$ locally uniformly in $\bar\Om\setminus\Om_*$. In turn, $\hat I=0$ on $\bar\Omega\setminus{\Om}_*$.

Suppose that  $\Om_*=\cup_{j=0}^m{\Om}_*^j$ with the properties as given by the assumptions.  It follows from \eqref{TH2.2-eq1} that $\int_{\mathring{\Om}_*}\hat{I}=N-|\Omega|r_{\min}>0$.

For each $1\leq j\leq m$,  by \eqref{Eq3-12} and the standard regularity theory of elliptic equations (e.g. \cite{GT}), restricting to a subsequence if necessary, we may assume that $I\to \hat{I}$ locally uniformly in $\mathring{\Om}_*^j$ as $d_I\to0$, where $\hat{I}\in C^2(\mathring{\Om}_*^j)\cap L^{\infty}(\Om_*^j)$. By \eqref{Eq3-9}-\eqref{Eq3-10}, we may further assume that
$$
\frac{d_S}{d_I}\Big(\frac{\kappa}{d_S}-r_{\min}\Big)\to \hat{a}_j\in(0,C]\ \ \mbox{as}\ d_I\to0.
$$
Moreover, $(\hat{a}_j,\hat{I})$ solves
\begin{equation}\label{Eq4-3}
        \Delta \hat{I}+\frac{\beta}{d_S}(\hat{a}_j-\hat{I})\hat{I}=0,\ \quad  x\in\mathring{\Om}_*^j.
\end{equation}
By the maximum principle, either $\hat I=0$ or $\hat I$ is a positive solution of \eqref{Eq4-3}. If $\hat I>0$,  noticing $\hat I=0$ on $\bar\Omega\setminus{\Om}_*$, $\hat I$ must satisfy \eqref{2-8aa}.
By $\int_{\mathring{\Om}_*}\hat{I}=N-|\Omega|r_{\min}$, \eqref{2-8ab} holds. This proves Proposition \ref{Prop2.1}.
\end{proof}



\subsection{Small  $d_S$:\ proof of Theorem \ref{TH2.3}(i)}
This subsection concerns the case  $d_S\to 0$ and proves Theorem \ref{TH2.3}(i).

\begin{proof}[Proof of Theorem \ref{TH2.3}(i)] Notice that  $\mathcal{R}_{0}$ is independent of $d_S$. Since $\mathcal{R}_{0}>1$,  by Theorem \ref{theorem_ex}, \eqref{1-1} admits at least one EE $(S,I)$ for any $d_S>0$. Let $({\kappa},\bar{S},\bar{I})$ be given as in Lemma \ref{Lem-2-1}.  Let
\begin{equation}\label{Eq9-2}
   C_*:= \limsup_{d_S\to 0}\frac{\kappa}{d_S}.
\end{equation}
By \eqref{Eq6-4}-\eqref{Eq6-5}, it is easily seen that $C_*\in[\frac{N}{|\Omega|},\infty]$.

 {\rm (i-a)} Suppose that $N<\int_{\Om}r^{\frac{1}{q}}$. We first claim that  $C_*$ is finite. Assume to the contrary that $C_*=\infty$. Then restricting to a subsequence if necessary, we have
\begin{equation}\label{Eq20-1}
    \lim_{d_S\to 0}\frac{\kappa}{d_{S}}=\infty.
\end{equation}
By  \eqref{Eq6-2}-\eqref{Eq6-3}, we have $S=\frac{\kappa}{d_{S}}(1-d_I\bar{I})$. Thus by \eqref{Eq6-5}, we find that
\begin{equation}\label{Eq9-1}
    \begin{cases}
\displaystyle        \frac{d_{S}}{\kappa}\Delta S +\beta( r-S^q)\bar{I}=0,\ \ \ & x\in\Om,\cr
        \p_{\nu} S=0, & x\in\p \Om.
    \end{cases}
\end{equation}
By Lemma \ref{l3.1}, we have
    $$
    \|S\|_{L^\infty(\Omega)}\le r_{\max}^{\frac{1}{q}},\ \ \ \quad \forall\; d_S> 0.
    $$
This and \eqref{Eq20-1} imply that
    $$
    \|1-d_I\bar{I}\|_{L^\infty(\Omega)}=\frac{d_{S}}{\kappa}\|S\|_{\infty}\le \frac{d_{S}}{\kappa}r_{\max}^{\frac{1}{q}}\to 0 \ \ \text{as}\ d_S\to 0.
    $$
Now one can employ the standard singular perturbation theory to  \eqref{Eq9-1} to deduce that
$S\to r^{\frac{1}{q}}$  uniformly on $\bar\Om$ as $d_S\to 0$. As a result,
\begin{equation}\label{Eq9-3}
\lim_{d_S\to 0}\int_{\Om}S=\int_{\Om}r^{\frac{1}{q}}>N,
\end{equation}
which contradicts the fact that $N>\int_{\Om}S$. Hence $C_*<\infty$.

By \eqref{Eq6-1} and $C_*<\infty$,  there is $C>0$ such that
$$
\|I\|_{L^\infty(\Omega)}\le \frac{1}{d_I}\kappa\le Cd_S,\quad \forall\, 0<d_S<1,
$$
which proves  \eqref{Eq8-03-1}.

Next, we show \eqref{Eq8-03-02}. Notice that $C_*<\infty$ and $0<\bar{I}<\frac{1}{d_I}$. Then by the elliptic estimates and \eqref{Eq6-5}, there exists $C>0$ such that $\|\bar I\|_{C^{2, \alpha}(\Omega)}\le C$ for all $0<d_S<1$. So up to a subsequence if necessary, $\frac{\kappa}{d_S}\to C^*$ and $\bar{I}\to I^*$ in $C^{1}(\bar{\Om})$ as $d_S\to0$, where $C^*\in[0,C_*]$ and  $I^*\in C^{2}(\Om)\cap C^{1}(\bar{\Om})$ with $0\le I^*\le \frac{1}{d_I}$. Moreover, $I^*$ is a classical solution of
\begin{equation}\label{3.29}
    \begin{cases}
        d_I\Delta I^* + \beta\Big(\big(C^*(1-d_II^*)\big)^q-r\Big)I^*=0,\ \ \ & x\in\Om,\cr
\mathfrak{}        \p_{\nu}I^*=0, & x\in\p\Om.
    \end{cases}
\end{equation}
Since $S=\frac{\kappa}{d_S}(1-d_I\bar{I})$, $S\to C^*(1-d_II^*)$  in $C^{1}(\bar{\Om})$ as $d_{S}\to 0$. Then by  \eqref{Eq8-03-1} and $\int_\Omega(S+I)=N$, $\int_{\Om}C^*(1-d_II^*)=N$. Hence, we have
 \begin{equation}\label{Eq8-23-0}
\frac{\kappa}{d_S}\to C^*=\frac{N}{\int_{\Om}(1-d_II^*)}.
 \end{equation}
 This shows \eqref{Eq8-03-02}.

 By \eqref{3.29} and the maximum principle, we have $I^*<\frac{1}{d_I}$ on $\bar{\Om}$. To complete the proof of (i-a), it remains to show that $I^*>0$. To see this,  let $\varphi$ be a positive eigenfunction corresponding to $\mathcal{R}_{0}$ which satisfies (see \cite{allen2008asymptotic})
 \begin{equation}\label{Eq8-23}
     \begin{cases}
\displaystyle         d_I\Delta\varphi +\beta\Big(\frac{1}{\mathcal{R}_{0}}\frac{N^q}{|\Om|^q}-r\Big)\varphi=0,\ \ \ & x\in\Om,\cr
         \p_{\nu}\varphi=0, & x\in\p\Om.
     \end{cases}
 \end{equation}
Define $\tilde{u}:=\frac{\bar{I}}{\|\bar{I}\|_{L^\infty(\Omega)}}$. Then it satisfies $0<\tilde{u}\le 1$, $\tilde{u}_{\max}=1$, and
 $$
  \begin{cases}
 \displaystyle   d_I\Delta \tilde{u} +\beta\Big(\Big(\frac{\kappa}{d_S}(1-d_I\bar{I})\Big)^q-r\Big)\tilde{u}=0,\ \ \  & x\in\Omega,\cr
    \p_{\nu}\tilde{u}=0, & x\in\partial\Omega.
    \end{cases}
 $$
By the  elliptic estimates, restricting to a subsequence if necessary, we have $\tilde{u}\to u^*$ as $d_S\to 0$ in $C^{1}(\bar{\Om})$, where $u^*$ is  a positive solution of
 $$
  \begin{cases}
    d_I\Delta {u}^* +\beta\Big( \Big(\frac{N(1-d_I{I^*})}{\int_{\Om}(1-d_II^*)}\Big)^q-r\Big){u}^*=0,\ \ \ & x\in\Omega,\cr
    \p_{\nu}{u}^*=0, & x\in\partial\Omega, \cr
    0<u^*\le 1,\ \ u^*_{\max}=1.
    \end{cases}
 $$
 Multiplying this equation by $\varphi$ and integrating over $\Omega$, we obtain
 \begin{equation*}
     d_I\int_{\Om}\nabla u^*\cdot\nabla\varphi+\int_{\Om}\gamma u^*\varphi=\int_{\Om}\beta \left[\frac{N(1-d_I{I^*})}{\int_{\Om}(1-d_II^*)}\right]^qu^*\varphi.
 \end{equation*}
 Similarly, multiplying \eqref{Eq8-23} by $u^*$ and integrating over $\Omega$, we can infer that
  \begin{equation*}
     d_I\int_{\Om}\nabla u^*\cdot\nabla\varphi+\int_{\Om}\gamma u^*\varphi=\int_{\Om}\frac{\beta}{\mathcal{R}_{0}} \left[\frac{N}{|\Om|}\right]^qu^*\varphi.
 \end{equation*}
 Therefore,
 $$\int_{\Om}\beta \left[\frac{N(1-d_I{I^*})}{\int_{\Om}(1-d_II^*)}\right]^qu^*\varphi= \int_{\Om}\frac{\beta}{\mathcal{R}_{0}} \left[\frac{N}{|\Om|}\right]^qu^*\varphi. $$
Since $u^*\varphi>0$ on $\bar{\Om}$, we have
$$
1-d_I\|I^*\|_{L^\infty(\Omega)}\le \frac{\int_{\Om}(1-d_II^*)}{\mathcal{R}_{0}^{\frac{1}{q}}|\Omega|}.
$$
 It follows that
$$
\|I^*\|_{L^\infty(\Omega)}\ge \frac{1}{d_{I}}\left(1-\frac{\int_{\Om}(1-d_II^*)}{\mathcal{R}_{0}^{\frac{1}{q}}|\Omega|}\right)
=\frac{\mathcal{R}_{0}^{\frac{1}{q}}-1}{d_I\mathcal{R}_{0}}^{\frac{1}{q}}
+\frac{\int_{\Om}I^*}{\mathcal{R}_{0}^{\frac{1}{q}}|\Om|}\ge \frac{\mathcal{R}_{0}^{\frac{1}{q}}-1}{d_I\mathcal{R}_{0}}^{\frac{1}{q}}>0.
$$
 This shows that $I^*\ge,\ne 0$, and hence $ \tilde{u}^*:=\beta(C^*(1-d_II^*))^qI^*\ge,\ne 0 $.  Note that the  linear operator $\{ u\in\cap_{p>1}W^{2,p}(\Omega) : \Delta u\in C(\bar{\Omega}) \ \text{and}\ \partial_{\nu}u=0 \ \text{on}\ \partial\Omega\}\ni u  \to \gamma u-d_I\Delta u\in C(\bar{\Omega})$ is invertible  since $ \gamma>0$ on $\bar{\Omega}$. Moreover, by the strong maximum principle and Hopf boundary lemma, its inverse function $(\gamma-d_I\Delta)^{-1}$ is strongly positive. Therefore, by \eqref{3.29}, $I^*=(\gamma-d_I\Delta)^{-1}\tilde{u}^*>0$ on $\bar{\Omega}$. 

{\rm (i-b)} Suppose that $N=\int_{\Om}r^{\frac{1}{q}}$. We distinguish two cases:  $C_*=\infty$ and $C_*<\infty$.

\textit{Case 1}. Suppose that  $C_*=\infty$. Then restricting to a subsequence if necessary, we may assume that \eqref{Eq20-1} holds. By the same arguments leading to \eqref{Eq9-3}, we get
$$
     \lim_{d_S\to0}\|S-r^{\frac{1}{q}}\|_{L^\infty(\Omega)}=0 \quad \text{and}\quad\lim_{d_S\to 0}\int_{\Om}S=\int_{\Om}r^{\frac{1}{q}}=N.
$$
Hence by $\int_\Omega (S+I)=N$,
 \begin{equation}\label{Eq11-2}
     \lim_{d_S\to 0}\int_{\Om}I=0.
 \end{equation}
 Observing that
 $$
 \sup_{d_S>0}\|\beta S^q-\gamma\|_{L^\infty(\Omega)}\le \|\beta\|_{L^\infty(\Omega)}\|r\|_{L^\infty(\Omega)}+\|\gamma\|_{L^\infty(\Omega)}<\infty,
 $$
 it follows from the second equation of \eqref{1-1} and the Harnack's inequality (see, e.g. \cite[Theorem 2.5]{huska2006harnack} or \cite[Lemma 2.2]{peng2008stationary}) that
 $$
 I_{\max}\le c_0 I_{\min},\ \ \ \quad \forall\; d_S>0
 $$
for some constant $c_0>0$. This together with \eqref{Eq11-2} implies that
$\|I\|_{L^\infty(\Omega)}\to 0$ as $d_S\to0$.

\textit{Case 2}. Suppose that $C_*<\infty$.  Then we are in the scenario of {\rm (i-a)}. Therefore, there is a  constant $C>0$ such that \eqref{Eq8-03-1} holds, which yields that $\|I\|_{L^\infty(\Omega)}\to 0$ as $d_S\to0$. Furthermore, \eqref{Eq8-03-02} holds up to a subsequence of $d_S\to 0$.

Hence in either \textit{Case 1} or \textit{Case 2}, we have  $\|I\|_{L^\infty(\Omega)}\to 0$ as $d_S\to0$. Moreover, up to a subsequence of $d_S\to 0$, either $S\to r^{\frac{1}{q}}$ or $S\to S^*$ uniformly on $\bar\Om$.

(i-c) Suppose that $N>\int_{\Om}r^{\frac{1}{q}}$. For each $n\ge 1$, consider the elliptic equation
 \begin{equation}\label{Eq12-4}
     \begin{cases}
         d_I\Delta w_n +\beta\Big[\Big(n+\frac{N}{|\Om|}\Big)^q(1-d_Iw_n)^q-r \Big]w_n=0,\ \ \ & x\in\Om,\cr
         \p_{\nu}w_n=0, & x\in\p\Om,\cr
         0<w_n<\frac{1}{d_I}.
     \end{cases}
 \end{equation}
Notice that $\mathcal{R}_{0}>1$  and the function
$$\big[0,\frac{1}{d_I}\big]\ni w\mapsto \Big(n+\frac{N}{|\Om|}\Big)^q(1-d_Iw)^q-r $$
is strictly decreasing. Then it follows from the classical results on the Fisher-KPP type equations that problem \eqref{Eq12-4} admits a unique positive solution, denoted by $w_n$, for every $n\ge 1$.  It is easy to check that the function
\begin{equation}\label{Eq12-1}
    S_{n}=\Big(n+\frac{N}{|\Om|}\Big)(1-d_Iw_n)
\end{equation}
satisfies
\begin{equation}\label{Eq12-2}
\begin{cases}
    \frac{1}{n+\frac{N}{|\Omega|}}\Delta S_n +\beta(r-S_n^q)w_n=0,\ \ \ & x\in\Om,\cr
    \p_{\nu}S_n=0, & x\in\p \Om.
\end{cases}
\end{equation}
By Lemma \ref{l3.1}, we have $\|S_n\|_{L^\infty(\Omega)}\le r_{\max}^{\frac{1}{q}},\ \forall n\ge 1$. This implies that
\begin{equation}\label{Eq12-6}
\|d_Iw_n-1\|_{L^\infty(\Omega)}=\frac{1}{n+\frac{N}{|\Om|}}\|S_n\|_{L^\infty(\Omega)}\le \frac{r_{\max}^{\frac{1}{q}}}{n+\frac{N}{|\Om|}}\to 0\ \ \text{as} \ n\to\infty.
\end{equation}
Hence, we can employ the singular perturbation theory to \eqref{Eq12-2} to deduce that
\begin{equation}\label{Eq12-7}
\|S_n-r^{\frac{1}{q}}\|_{L^\infty(\Omega)}\to 0\ \ \text{as}\ n\to \infty.
\end{equation}

Now for all $n\ge 1$, let
\begin{equation}\label{SD1}
d_{S,n}=\frac{N-\int_{\Om}S_n}{\Big(n+\frac{N}{|\Om|}\Big)\int_{\Om}w_n} \quad \text{and}\quad I_n=\frac{(N-\int_{\Om}S_n)}{\int_{\Om}w_n}w_n.
\end{equation}
It follows from \eqref{Eq12-6} and \eqref{Eq12-7} that $d_{S,n}\to 0$ and $\|I_n-{(N-\int_{\Om}r^{\frac{1}{q}})/|\Omega|}\|_{L^\infty(\Omega)}\to 0$ as $n\to\infty$. Note also by \eqref{Eq12-4} and \eqref{Eq12-2} that $(S_n,I_n)$ is an EE of \eqref{1-1} with $d_S=d_{S,n}$ for all $n\ge 1$.
\end{proof}

\begin{rk}\label{r-th2.3-a} Suppose that the hypotheses of Theorem \ref{TH2.3}(i) hold and fix $d_I>0$ such that $\mathcal{R}_{0}>1$.
\begin{itemize}
    \item[\rm(1)] If $N=\int_{\Om}r^{\frac{1}{q}}$, Theorem \ref{TH2.3}{(i-b)} shows that for any EE $(S,I)$ of \eqref{1-1}, either $S\to S^*$ or $ S\to r^{\frac{1}{q}}$ uniformly on $\bar\Om$ up to a subsequence as $d_S\to 0$, where $S^*$ is given as in \eqref{Eq8-03-02}. It is unclear which one of these scenarios holds. However, if we further assume that $\beta$ is constant and $q=1$, then we can find a sequence of solutions such that the later scenario occurs. Indeed, in this case, consider the sequence $\{w_n\}_{n\ge 1}$ of positive solutions of \eqref{Eq12-4} constructed in the proof of Theorem \ref{TH2.3}{(i-b)}. Dividing the first equation in \eqref{Eq12-4} by $\beta w_n$ and integrating the resulting equation gives
    \begin{align}
    0= d_I\int_{\Om}\frac{|\nabla \ln w_n|^2}{\beta}+\int_{\Om}S_n-\int_{\Om}r=\frac{d_I}{\beta}\int_{\Om}|\nabla \ln w_n|^2-\Big(N-\int_{\Om}S_n\Big),
    \end{align}
    where $S_n$ is defined as in \eqref{Eq12-1}. This shows that $d_{S,n}$ defined in \eqref{SD1} is positive for every $n\ge 1$. Therefore, $(S_n,I_n)$ defined as in \eqref{Eq12-1} and \eqref{SD1} is a sequence of EEs of \eqref{1-1} with $d_S=d_{S,n}$ satisfying $S_n\to r^{\frac{1}{q}}$ and $I_n\to 0$ uniformly on $\bar\Om$ as $n\to\infty$.

    \item[\rm(2)] To show that the limit in Theorem \ref{TH2.3}(i-c) holds for any EE (at least under extra assumptions), it is not hard to see that we just need to rule out $C_*<\infty$ in the proofs of Theorem \ref{TH2.3}(i-a) and (i-b), which can be achieved by finding conditions to guarantee that \eqref{Eq8-03-03} has no positive solution. Note that if $I^*$ is a positive solution of \eqref{Eq8-03-03} and $x_1\in\bar{\Omega}$ such that $I^*(x_1)=I^*_{\min}$, then by the maximum principle (Lemma \ref{l3.1}), it holds that
$$
N^q\frac{(1-d_II^*(x_1))^q}{\|1-d_II^*\|_{L^1(\Omega)}^q}\le r(x_1) \le r_{\max},
$$
which leads to
$$
1-d_II^*(x)\le 1-d_II_{\min}^*\le \|1-d_II^*\|_{L^1(\Omega)}\frac{r^{\frac{1}{q}}_{\max}}{N},\quad \forall\, x\in\bar\Omega.
$$
This implies that
\begin{equation}\label{JH1}
\|1-d_II^*\|_{L^1(\Omega)}\le\|1-d_II_{\min}^*\|_{L^1(\Omega)}\leq \|1-d_II^*\|_{L^1(\Omega)}|\Om|\frac{r^{\frac{1}{q}}_{\max}}{N},
\end{equation}
which is impossible if ${N}/{|\Omega|}>r_{\max}^{\frac{1}{q}}$ since $d_II^*<1$.
If ${N}/{|\Omega|}=r_{\max}^{\frac{1}{q}}$, by  \eqref{JH1} and $0<d_II^*<1 $ again, we have  $I^*\equiv I^*_{\min}$. In view of \eqref{Eq8-03-03}, $r$ is  constant. This is also impossible since $1<{\mathcal{R}_{0}}\le {\big(\frac{N}{|\Omega|}\big)^{q}}/{r_{\min}}={r_{\max}}/{r_{\min}}=1$.  Therefore if $\frac{N}{|\Omega|}\ge r_{\max}^{q}$, then $(S,I)\to (r^{\frac{1}{q}},\ \frac{1}{|\Om|}(N-\int_{\Om}r^{\frac{1}{q}}))$ uniformly on $\bar\Om$ as $d_S\to 0$.

\end{itemize}
\end{rk}

\subsection{Small $d_S$ and $d_I$:\ proof of Theorem \ref{TH2.4}{\rm (i)}}
This subsection is devoted to the case that $d_S\to 0$ and $d_I/d_S\to \sigma$ and  the proof of Theorem \ref{TH2.4}(i).

\begin{proof}[Proof of Theorem \ref{TH2.4}{\rm (i)}] 
The existence of EE can be guaranteed by the same argument of Theorem \ref{TH2.1}{\rm (i)}. Let $(\kappa,\bar{S},\bar{I})$ be given by Lemma \ref{Lem-2-1}(i).
Note from \eqref{Eq6-1} that
\begin{equation}\label{LU1}
\frac{\kappa}{d_S}=\frac{1}{|\Omega|}\int_{\Omega}\Big(S+\frac{d_I}{d_S}I\Big)\le \frac{N}{|\Omega|}\Big(1+\frac{d_I}{d_S}\Big).
\end{equation}
Let $\tilde{u}=d_S\bar{I}$. By \eqref{Eq6-5}, $\tilde{u}$ solves
\begin{equation}\label{LU2}
    \begin{cases}
\displaystyle d_I\Delta\tilde{u}+\beta\Big[\Big(\frac{\kappa}{d_S}\Big)^{q}\Big(1-\frac{d_I}{d_S}\tilde{u}\Big)^q-r\Big]\tilde{u}=0,\ \ \ & x\in\Om,\cr
\p_{\nu}\tilde{u}=0, & x\in\p\Om
    \end{cases}
\end{equation}
and
\begin{equation}\label{LU3}
    \tilde{u}=\frac{d_S}{d_I}d_I\bar{I}<\frac{d_S}{d_I}.
\end{equation}

We claim that
\begin{equation}\label{LU4}
    \lim_{d_I+\big|\frac{d_I}{d_S}-\sigma\big|\to0}
    \Big(\Big|\frac{\kappa}{d_S}-\tilde{\kappa}_{\sigma}\Big|+\Big\|\tilde{u}
    -\Big(1-\frac{r^{\frac{1}{q}}}{\tilde{\kappa}_{\sigma}}\Big)_+\Big\|_{L^\infty(\Omega)}\Big)=0,
\end{equation}
where the constant $\tilde{\kappa}_{\sigma}$ is the unique positive solution of the algebraic equation \eqref{TH2.4-eq2}. Indeed, by \eqref{Eq8-2} and \eqref{LU1}, there is ${\kappa}^*\in\big[r_{\min}^{\frac{1}{q}},\frac{N(1+\sigma)}{|\Omega|}\big]$ such that, restricted  to a subsequence if necessary, we have
$$
\frac{\kappa}{d_S}\to {\kappa}^*\quad \text{as}\  d_I+|\frac{d_I}{d_S}-\sigma|\to 0.
$$
 Moreover, due to \eqref{LU3}, it holds that
 $$
 \sup_{d_I+|\frac{d_I}{d_S}-\sigma|\le \frac{\sigma}{2}}\|\tilde{u}\|_{L^\infty(\Omega)}\le 1+\frac{2}{\sigma}.
 $$
So we can employ the singular perturbation theory for elliptic equations and deduce from \eqref{LU2} that
$$
\mbox{$\tilde{u}\to \frac{1}{\sigma}\Big(1-\frac{r^{\frac{1}{q}}}{{\kappa}^*}\Big)_{+}$\ \ \ \ uniformly on\ $\bar\Om$\,\ \ as\ $d_I+|\frac{d_I}{d_S}-\sigma|\to 0$.}
$$
By \eqref{Eq6-4}, we have
$$
N=\int_{\Om}\left[\frac{\kappa}{d_S}\left(1-\frac{d_I}{d_S}\tilde{u}\right)+\frac{\kappa}{d_S}\tilde{u}\right].
$$
Taking $d_I+|\frac{d_I}{d_S}-\sigma|\to0$, we see that ${\kappa}^*$ solves the algebraic equation \eqref{TH2.4-eq2}. Since \eqref{TH2.4-eq2} has a unique positive solution $\tilde{\kappa}_{\sigma}$,  we have  ${\kappa}^*=\tilde{\kappa}_{\sigma}$. Since $\tilde{\kappa}_{\sigma}$ is independent of the chosen subsequence, \eqref{LU4} holds.

By \eqref{LU4}, we have
\begin{equation}\label{LU4-aa}
S=\frac{\kappa}{d_S}(1-\frac{d_I}{d_S}\tilde{u})\to \tilde{\kappa}_{\sigma}\Big[1-\Big(1-\frac{r^{\frac{1}{q}}}{\tilde{\kappa}_{\sigma}}\Big)_+\Big]
=\min\{\tilde{\kappa}_{\sigma},r^{\frac{1}{q}}\}=S_{\sigma}
\end{equation}
and
\begin{equation}\label{LU4-ab}
I=\frac{\kappa}{d_S}\tilde{u}\to \frac{\tilde{\kappa}_{\sigma}}{\sigma}\Big(1-\frac{r^{\frac{1}{q}}}{\tilde{\kappa}_{\sigma}}\Big)_+
=\frac{1}{\sigma}\big(\tilde{\kappa}_{\sigma}-r^{\frac{1}{q}}\big)_+=I_{\sigma}
\end{equation}
uniformly on $\bar\Om$ as $d_I+|\frac{d_I}{d_S}-\sigma|\to 0$. Hence \eqref{TH2.4-eq1} holds.

To complete the proof of the theorem, it remains to show that $(\tilde{\kappa}_{\sigma},S_{\sigma},I_{\sigma})$ satisfies {\rm(i-a)--(i-d)}.

\quad {\rm (i-a)} Suppose that $N<\int_{\Omega}r^{\frac{1}{q}}$. It follows from \eqref{TH2.4-eq2}, \eqref{Eq8-2} and $\frac{N}{|\Om|}>r_{\min}^{\frac{1}{q}}$ that $\tilde{\kappa}_{\sigma}\in (r_{\min}^{\frac{1}{q}},r_{\max}^{\frac{1}{q}})$ for all $\sigma>0$. Hence
$$
\frac{1}{\sigma}\int_{\Om}(\tilde{\kappa}_{\sigma}-r^{\frac{1}{q}})_+\to0 \ \ \mbox{as}\ \sigma\to\infty.
$$
Moreover, restricted to a subsequence if necessary, we have  $\tilde{\kappa}_{\sigma}\to {\kappa}^*\in[r_{\min}^{\frac{1}{q}},r^{\frac{1}{q}}_{\max}]$ as $\sigma\to\infty$. By \eqref{TH2.4-eq2}, $N=\int_{\Om}\min\{{\kappa}^*,r^{\frac{1}{q}}\}$. Since $N<\int_{\Om}r^{\frac{1}{q}}$ and $\frac{N}{|\Om|}>r_{\min}^{\frac{1}{q}}$, there is a unique $\tilde{\kappa}_{\infty}\in(r_{\min}^{\frac{1}{q}},r^{\frac{1}{q}}_{\max})$ satisfying $N=\int_{\Om}\min\{\tilde{\kappa}_{\infty},r^{\frac{1}{q}}\}$. Therefore, we have that $\tilde{\kappa}_{\sigma}\to \tilde{\kappa}_{\infty}$ as $\sigma\to\infty$. Furthermore, it follows from \eqref{LU4-aa}-\eqref{LU4-ab} that
$$
\mbox{$S_{\sigma}\to \min\{\tilde{\kappa}_{\infty},r^{\frac{1}{q}}\}$\ \, and\ $I_{\sigma}\to 0$\ \ \ \ uniformly on\ $\bar\Om$\ \ as\ $\sigma\to\infty$.}
$$

\quad {\rm (i-b)} Suppose that $N=\int_{\Om}r^{\frac{1}{q}}$. As in (i-a), we have $\tilde{\kappa}_{\sigma}\in (r_{\min}^{\frac{1}{q}},r_{\max}^{\frac{1}{q}})$ for all $\sigma>0$. Hence, we can proceed as in {\rm(i-a)} to show that $\tilde{\kappa}_{\sigma}\to r_{\max}^{\frac{1}{q}}$, $I_{\sigma}\to 0$ and $S_{\sigma}\to r^{\frac{1}{q}}$ uniformly on $\bar\Om$ as $\sigma\to\infty$.

\quad {\rm (i-c)} Suppose that $N>\int_{\Omega}r^{\frac{1}{q}}$. By \eqref{TH2.4-eq2},
$$
\tilde{\kappa}_{\sigma}>\frac{1}{|\Omega|}\int_{\Omega}\Big(\tilde{\kappa}_{\sigma}-r^{\frac{1}{q}}\Big)_+\ge \frac{\sigma}{|\Omega|}\Big(N-\int_{\Omega}r^{\frac{1}{q}}\Big)\to \infty\ \ \,\text{as} \ \sigma\to\infty.
$$
By \eqref{TH2.4-eq2} again, there is $\sigma^*\gg 1$ such that
\begin{equation}\label{LU4-ac}
\tilde{\kappa}_{\sigma}=\frac{\sigma}{|\Omega|}\left(N-\int_{\Om}r^{\frac{1}{q}}\right)
+\frac{1}{|\Omega|}\int_{\Omega}r^{\frac{1}{q}}\ \,\ \ \ \mbox{whenever}\ \sigma\ge \sigma^*.
\end{equation}
 It follows from \eqref{LU4-aa}-\eqref{LU4-ab} that $S_{\sigma}\to r^{\frac{1}{q}} $ and $I_{\sigma}\to \frac{1}{|\Omega|}\Big(N-\int_{\Om}r^{\frac{1}{q}}\Big)$ uniformly on $\bar\Om$ as $\sigma\to\infty$.

\quad {\rm (i-d)} By \eqref{TH2.4-eq2}, $\int_{\Om}\big(\tilde{\kappa}_{\sigma}-r^{\frac{1}{q}}\big)_+\le \sigma N\to 0$ as $\sigma\to0$. Hence,  it follows from \eqref{Eq8-2} that $\tilde{\kappa}_{\sigma}\to r_{\min}^{\frac{1}{q}}$ as $\sigma\to0$. By \eqref{LU4-aa}-\eqref{LU4-ab}, $S_{\sigma}\to r_{\min}^{\frac{1}{q}}$ uniformly on $\bar\Om$ and $I_{\sigma}\to 0$ uniformly on $K$ for every compact subset $K\subset \bar\Om\setminus\Om_*$, and $\int_{\Om}I_{\sigma}\to N-|\Omega|r^{\frac{1}{q}}_{\min}$ as $\sigma\to0$.

Suppose in addition that $|\Om_*|>0$. Rewriting \eqref{TH2.4-eq2} as

$$
N=r_{\min}^{\frac{1}{q}}|\Om_*|+\int_{\Om\setminus\Om_*}\min\{\tilde{\kappa}_{\sigma},r^{\frac{1}{q}}\}+\frac{(\tilde{\kappa}_{\sigma}-r_{\min}^{\frac{1}{q}})}{\sigma}|\Om_*|+\int_{\Omega\setminus\Om_*} \frac{(\tilde{\kappa}_{\sigma}-r^{\frac{1}{q}})_+}{\sigma},
$$
we obtain that
\begin{equation}\label{LU8}
  \frac{(\tilde{\kappa}_{\sigma}-r_{\min}^{\frac{1}{q}})}{\sigma}=\frac{1}{|\Om_*|}\left( N-r_{\min}^{\frac{1}{q}}|\Om_*|-\int_{\Om\setminus\Om_*}\min\{\tilde{\kappa}_{\sigma},r^{\frac{1}{q}}\}-\int_{\Omega\setminus\Om_*} \frac{(\tilde{\kappa}_{\sigma}-r^{\frac{1}{q}})_+}{\sigma}\right).
\end{equation}
Hence, for every $x\in\Om_*$,
\begin{eqnarray*}
    I_{\sigma}(x)&=&\frac{1}{\sigma}(\tilde{\kappa}_{\sigma}-r^{\frac{1}{q}}_{\min})\\
    &=&\frac{1}{|\Om_*|}\left( N-r_{\min}^{\frac{1}{q}}|\Om_*|-\int_{\Om\setminus\Om_*}\min\{\tilde{\kappa}_{\sigma},r^{\frac{1}{q}}\}-\int_{\Omega\setminus\Om_*} \frac{(\tilde{\kappa}_{\sigma}-r^{\frac{1}{q}})_+}{\sigma}\right).
\end{eqnarray*}
Therefore, it holds that
\begin{eqnarray}\ \ \ \label{LU9}
    &&\Big\|I_{\sigma}-\frac{(N-r_{\min}^{\frac{1}{q}}|\Om|)}{|\Om_*|}\Big\|_{L^{\infty}(\Om_*)}\le \\
    &&\hspace{1.5cm}\frac{1}{|\Om_*|}\int_{\Om\setminus\Om_*}
    \min\big\{\tilde{\kappa}_{\sigma}-r_{\min}^{\frac{1}{q}},r^{\frac{1}{q}}-r_{\min}^{\frac{1}{q}}\big\}
    +\frac{1}{|\Om_*|}\int_{\Omega\setminus\Om_*} \frac{(\tilde{\kappa}_{\sigma}-r^{\frac{1}{q}})_+}{\sigma}.\nonumber
\end{eqnarray}
Note from \eqref{LU8} that
\begin{eqnarray*}
\int_{\Om\setminus\Om_*}\min\{\tilde{\kappa}_{\sigma}-r_{\min}^{\frac{1}{q}},r^{\frac{1}{q}}-r_{\min}^{\frac{1}{q}}\}&\le&\int_{\Om\setminus\Om_*}(\tilde{\kappa}_{\sigma}-r_{\min}^{\frac{1}{q}})\\
&\le& \frac{\sigma}{|\Om_*|}(N-r_{\min}^{\frac{1}{q}})\to 0\ \ \text{as}\ \sigma\to0,
\end{eqnarray*}
and
$$
0\le \frac{(\tilde{\kappa}_{\sigma}-r^{\frac{1}{q}})_+}{\sigma}\le \frac{(\tilde{\kappa}_{\sigma}-r^{\frac{1}{q}}_{\min})}{\sigma}\le \frac{1}{|\Om_*|}(N-r_{\min}^{\frac{1}{q}}).
$$
Since $r_{\min}^{\frac{1}{q}}<r^{\frac{1}{q}}(x)$ for every $x\in\Om\setminus\Om_*$ and $\tilde{\kappa}_{\sigma}\to r_{\min}^{\frac{1}{q}}$ as  $\sigma\to0$,
we have
$$
\frac{(\tilde{\kappa}_{\sigma}-r^{\frac{1}{q}}(x))_+}{\sigma}\to 0, \quad \forall x\in\Om\setminus\Om_*\ \ \text{as}\ \sigma\to 0.
$$
By the Lebesgue dominated convergence theorem,
$$
\int_{\Omega\setminus\Om_*} \frac{(\tilde{\kappa}_{\sigma}-r^{\frac{1}{q}})_+}{\sigma}\to 0\ \ \mbox{as}\ \sigma\to0.
$$
As a result, we get from \eqref{LU9} that
$$ \Big\|I_{\sigma}-\frac{(N-r_{\min}^{\frac{1}{q}}|\Om|)}{|\Om_*|}\Big\|_{L^{\infty}(\Om_*)} \to 0\ \ \mbox{ as}\ \sigma\to 0.
$$
The proof is now complete.
\end{proof}

\section{Proofs of the main results: case of $0<p<1$}\label{case p<1} In this section, we assume that $0<p<1$ and prove Theorem \ref{TH2.1}(ii), Theorem \ref{TH2.3}(ii) and Theorem \ref{TH2.4}(ii). By Theorem \ref{theorem_ex},  \eqref{1-1} has an EE $(S, I)$ for all $d_S, d_I>0$. In this section, let $({\kappa},\bar{S},\bar{I})$ be given as in Lemma \ref{Lem-2-1}.

\subsection{Small $d_I$: proof of Theorem \ref{TH2.1}{\rm (ii)}}
We need the following two lemmas.

\begin{lem}\label{Lem-2-3} Suppose that $0<p<1$ and $q>0$. Let $(S,I)$ be an EE of \eqref{1-1}. Then it holds that
\begin{equation}\label{Eq6-6}
I_{\max}\leq \Big[\Big(\frac{\beta}{\gamma}\Big)_{\max}S_{\max}^q\Big]^{\frac{1}{1-p}}
\end{equation}
and
\begin{equation}\label{Eq6-6s}
S_{\max}\leq \Big[r_{\max}I_{\max}^{1-p}\Big]^{\frac{1}{q}}.
\end{equation}

\end{lem}
\begin{proof} Set $I(x_0)=I_{\max}$ for some $x_0\in\bar\Omega$. Then by Lemma \ref{l3.1}, we have $\Delta I(x_0)\leq0$. Thus, the second equation in \eqref{1-1} indicates that
$$
\beta(x_0) S^q(x_0)I^p(x_0)-\gamma(x_0)I(x_0)\geq0,
$$
from which it follows that \eqref{Eq6-6} holds.

Similarly, suppose that $S(y_0)=S_{\max}$ for some $y_0\in\bar\Omega$. Then Lemma \ref{l3.1} asserts
$$
-\beta(y_0)S^q(y_0)I^p(y_0)+\gamma(y_0)I(y_0)\ge 0.
$$
This implies \eqref{Eq6-6s}.
\end{proof}

\begin{lem}\label{Lem-2-4} Suppose that $0<p<1$ and $q>0$.  Let $(S,I)$ be an EE of \eqref{1-1} and $\kappa$ be given as in Lemma \ref{Lem-2-1}. Then we have
\begin{equation}\label{Eq6-7}
    \frac{\kappa}{d_S}\ge N_*,
\end{equation}
where $N_*$ is the unique positive solution of the algebraic equation
$$  \frac{N}{|\Om|}= N_*+\Big(\frac{\beta}{\gamma}\Big)_{\max}^{\frac{1}{1-p}}N_*^{\frac{q}{1-p}}.$$
Furthermore, if $d_S\ge d_I$, then we have
\begin{equation}\label{Eq6-8}
    \frac{\kappa}{d_S}\le \frac{N}{|\Omega|}.
\end{equation}
\end{lem}
\begin{proof}
It follows from \eqref{Eq6-1} that $S<\frac{\kappa}{d_S}$. Hence
    \begin{align*}
        N=\int_{\Om}S+\int_{\Om}I\le \left(\frac{\kappa}{d_S}+I_{\max}\right)|\Omega|.
    \end{align*}
This together with \eqref{Eq6-6} yields
    $$
    \frac{N}{|\Om|}\le \frac{\kappa}{d_S}+\Big(\frac{\beta}{\gamma}\Big)_{\max}^{\frac{1}{1-p}}\Big(\frac{\kappa}{d_S}\Big)^{\frac{q}{1-p}}.
    $$
Since the function
    $$[0,\infty)\ni s\mapsto s+\Big(\frac{\beta}{\gamma}\Big)_{\max}^{\frac{1}{1-p}}s^\frac{q}{1-p}
    $$ is strictly increasing, continuous,
    $$ \lim_{s\to 0^+}\Big(s+\Big(\frac{\beta}{\gamma}\Big)_{\max}^{\frac{1}{1-p}}s^\frac{q}{1-p}\Big)=0 \quad \text{and}\quad \lim_{s\to\infty}\Big(s+\Big(\frac{\beta}{\gamma}\Big)_{\max}^{\frac{1}{1-p}}s^\frac{q}{1-p}\Big)=\infty,
    $$
     \eqref{Eq6-7} holds.
    It is clear that \eqref{Eq6-8} follows from \eqref{Eq6-4}.
\end{proof}

We are now ready to prove Theorem \ref{TH2.1}(ii).

\begin{proof}[Proof of Theorem \ref{TH2.1}{\rm (ii)}]
 By \eqref{Eq6-1}, we have $S< \frac{\kappa}{d_S}$. Then combining \eqref{Eq6-6} and \eqref{Eq6-8},
 \begin{equation}\label{Eq6-9}
    \|I\|_{L^\infty(\Omega)}\le C_1:=\Big[\Big(\frac{\beta}{\gamma}\Big)_{\max}\frac{N^q}{|\Om|^q}\Big]^{\frac{1}{1-p}}, \ \quad\ \forall\; 0<d_I\le d_S.
  \end{equation}
Since $\bar{I}=\frac{I}{\kappa}=\frac{d_S}{\kappa}\cdot\frac{I}{d_S}$,
we can obtain from \eqref{Eq6-7} and \eqref{Eq6-9} that
\begin{equation}\label{Eq6-10}
    \|\bar{I}\|_{L^\infty(\Omega)}\le \frac{C_1}{d_SN_*}, \ \quad \forall \;0<d_I\le d_S,
\end{equation}
where $N_*$ is the positive constant given in \eqref{Eq6-7}.

We now claim that
    \begin{equation}\label{Eq6-11}
        \lim_{d_I\to0}\|S-\frac{\kappa}{d_S}\|_{L^\infty(\Omega)}=0 \quad \text{and}\quad \lim_{d_I\to0}d_I\|\bar{I}\|_{L^\infty(\Omega)}=0.
    \end{equation}
Indeed, using \eqref{Eq6-1} and \eqref{Eq6-9}, we have
    $$
    \|\frac{\kappa}{d_S}-S\|_{L^\infty(\Omega)}=\frac{d_I}{d_S}\|I\|_{L^\infty(\Omega)}\le \frac{d_I}{d_S}C_1\to 0 \quad \text{as}\ d_I\to0.
    $$
Observe also from \eqref{Eq6-10} that
$$
d_I\|\bar{I}\|_{L^\infty(\Omega)}\le\frac{d_IC_1}{d_SN_*}\to 0\,\ \ \mbox{as}\ d_I\to 0.
$$
Thus, \eqref{Eq6-11} holds.

Restricting to a subsequence if necessary, we may assume $\frac{\kappa}{d_S}\to \tilde{S}_*$ for some $\tilde{S}_*\in [N_*,\frac{N}{|\Omega|}]$ as $d_I\to 0$.
Then we further claim the following:
\begin{equation}\label{Eq6-12}
\lim_{d_I\to0}\Big\|\bar{I}
-\Big(\frac{\beta}{\gamma}\Big)^{\frac{1}{1-p}}\Big(\frac{\tilde{S}_*^{q+p-1}}
{d_S^{1-p}}\Big)^{\frac{1}{1-p}}\Big\|_{L^\infty(\Omega)}=0
\end{equation}
and
\begin{equation}\label{Eq6-13}
N=\int_{\Om}\Big[\tilde{S}_* + \Big(\frac{\beta}{\gamma}\Big)^{\frac{1}{1-p}}\tilde{S}_*^{\frac{q}{1-p}} \Big].
\end{equation}
Indeed, using $\kappa^p\to d_S^p\tilde{S}_*^p>0$ as $d_I\to 0$ and  \eqref{Eq6-10},  we can employ the singular perturbation theory to deduce from \eqref{Eq6-5} that \eqref{Eq6-12} holds. Letting $d_I\to0$ in \eqref{Eq6-4}, we obtain \eqref{Eq6-13}.
Since this algebraic equation has a unique positive solution $S_*$, we have $\tilde S_*=S_*$. Since  $\tilde S_*$ is independent of the chosen subsequence, $\frac{\kappa}{d_S}\to S_*$ as $d_I\to 0$. By \eqref{Eq6-11},  $S\to S_*$ uniformly on $\bar{\Om}$ as $d_I\to0$. Moreover, thanks to \eqref{Eq6-12} we obtain that $$
I=d_S\frac{\kappa}{d_S}\bar{I}\to d_S S_*\Big(\frac{\beta}{\gamma}\frac{S_*^{q+p-1}}{d_S^{1-p}}\Big)^{\frac{1}{1-p}} =\Big(\frac{\beta}{\gamma}S_*^{q}\Big)^{\frac{1}{1-p}}
$$
uniformly on $\bar{\Om}$ as $d_I\to 0$.
\end{proof}

\subsection{Small $d_S$: proof of Theorem \ref{TH2.3}{\rm (ii)}}
The following result is similar to Lemma \ref{Lem-2-4}, so we omit the proof here.
\begin{lem}\label{Lem-2-4s} Suppose that $0<p<1$ and $q>0$.  Let $(S,I)$ be an EE of \eqref{1-1} and $\kappa$ be given as in Lemma \ref{Lem-2-1}. Then we have
\begin{equation}\label{Eq6-7s}
    \frac{\kappa}{d_I}\ge M_*,
\end{equation}
where $M_*$ is the unique positive solution of the algebraic equation
$$
\frac{N}{|\Om|}= M_*+r_{\max}^{\frac{1}{q}}M_*^{\frac{1-p}{q}}.
$$
Furthermore, if $d_S\le d_I$, then we have
\begin{equation}\label{Eq6-8s}
    \frac{\kappa}{d_I}\le \frac{N}{|\Omega|}.
\end{equation}
\end{lem}

With Lemma \ref{Lem-2-4s}, we are able to prove Theorem \ref{TH2.3}{\rm (ii)}.
\begin{proof}[Proof of Theorem \ref{TH2.3}{\rm (ii)}]  We first claim that $\lim_{d_S\to 0} \|\bar I-1/d_I\|_{L^\infty(\Omega)}=0$. To see this, by \eqref{Eq6-8s}, we have ${\kappa}\le \frac{N}{|\Omega|}d_I$ for all $0<d_S\le d_I$. Thus by the third equation of \eqref{Eq6-5},
\begin{equation}\label{Eq6-8s-1}
    I=\kappa\bar I\le \frac{N}{|\Omega|}, \quad \forall\;0<d_S\le d_I.
\end{equation}
Hence, by \eqref{Eq6-6s}, we obtain
\begin{equation}\label{Eq6-8s-2}
S\leq \Big[r_{\max}\Big(\frac{N}{|\Omega|}\Big)^{1-p}\Big]^{\frac{1}{q}},\ \ \ \ \forall\; 0<d_S\le d_I.
\end{equation}

By the definition of $\bar S$ and $\bar I$, we know that $d_S\bar S+d_I\bar I=1$, and so
$1-d_I\bar I=d_S\bar S=\frac{d_S}{\kappa}S$. It then follows that
\begin{eqnarray*}
\|1-d_I\bar I\|_{L^\infty(\Omega)}&\leq&\frac{d_S}{\kappa}\|S\|_{L^\infty(\Omega)}\\
&=&\frac{d_I}{\kappa}
\cdot\frac{d_S}{d_I}\|S\|_{L^\infty(\Omega)}\\
&\leq&\frac{d_S}{d_IM_*}\Big[r_{\max}\Big(\frac{N}{|\Omega|}\Big)^{1-p}\Big]^{\frac{1}{q}}\to0 \ \ \text{as}\ d_S\to 0,
\end{eqnarray*}
where we used \eqref{Eq6-7s} and \eqref{Eq6-8s-2}.
Consequently, together with \eqref{Eq6-8s}, we obtain
$$
\lim_{d_S\to 0} \|I-\kappa/d_I\|_{L^\infty(\Omega)}=\lim_{d_S\to 0} \kappa\|\bar I-1/d_I\|_{L^\infty(\Omega)}=0.
$$

By \eqref{Eq6-8s-1}-\eqref{Eq6-8s-2}, there exists $C>0$ such that $S,\ I<C$ on $\bar\Omega$ for all $0<d_S\le d_I$. Since $I$ satisfies
   $$
    \begin{cases}
    d_I\Delta I +\beta S^qI^p-\gamma I=0,\ \ \ & x\in\Om,\cr
    \p_{\nu}I=0, & x\in \p\Om,
    \end{cases}
    $$
the well-known elliptic estimates imply that $\|I\|_{W^{2, \ell}(\Omega)}<C$ by taking $C$ to be larger if necessary ($\ell>n$ is fixed) for all  $0<d_S\le d_I$. So by the standard embedding theorems, up to a subsequence if necessary,
$$
\lim_{d_S\to 0}\kappa= \kappa_* \ \ \text{and} \ \  \lim_{d_S\to 0} \|I-I_*\|_{L^\infty(\Omega)}=0
$$
for some $\kappa_*>0$ and $I_*=\kappa_*/d_I$ as $d_S\to 0$. Observing that  $S$ solves
   $$
    \begin{cases}
    d_{S}\Delta S -\beta S^qI^p+\gamma I=0, \ \ \ & x\in\Om,\cr
    \p_{\nu} S=0, & x\in \p\Om,
    \end{cases}
    $$
the standard singular perturbation techniques can be used to verify $\lim_{d_S\to 0} \|S-S_*\|_{L^\infty(\Omega)}=0$ with $S_*$ fulfilling
$-\beta S_*^qI_*^p+\gamma I_*=0$, and in turn $S_*=r^\frac{1}{q} I_*^\frac{1-p}{q}$. Due to $\int_\Omega (S+I)=N$, $I_*$ satisfies \eqref{Istr}.
Since \eqref{Istr} has a unique solution $I^*$, the limits $I^*$ and $\kappa_*$ are independent of the chosen subsequence. Hence, \eqref{SIc} holds.
\end{proof}

\subsection{Small $d_S$ and $d_I$: proof of Theorem \ref{TH2.4}(ii)}
The following result will be used in the proof of Theorem \ref{TH2.4}(ii).

\begin{lem}\label{lem-6} Let $0<p<1$, $q>0$ and  $\sigma>0$ be fixed. Then for every $\tilde{\kappa}\ge 0$ and $x\in\bar{\Om}$, there exists a unique $I_{\sigma}(x,\tilde{\kappa})\ge 0$ satisfying the algebraic equation
\begin{equation}\label{LU2-1}
    \tilde{\kappa}=\sigma I_{\sigma}(x,\tilde{\kappa}) +r^{\frac{1}{q}}(x)[I_{\sigma}(x,\tilde{\kappa})]^{\frac{1-p}{q}}.
\end{equation}
Moreover, the mapping $\bar{\Om}\times[0,\infty)\ni(x,\tilde{\kappa})\mapsto I_{\sigma}(x,\tilde{\kappa}) $ is continuous, strictly increasing in $\tilde{\kappa}$, $I_{\sigma}(x,0)\equiv 0$, and there is a unique positive constant $\tilde{\kappa}_{\sigma}$ satisfying the algebraic equation
\begin{equation}\label{LU2-2}
    N=\tilde{\kappa}_{\sigma}|\Om|+(1-\sigma)\int_{\Om}I_{\sigma}(x,\tilde{\kappa}_{\sigma}).
\end{equation}

\end{lem}
\begin{proof} Observe that the function
$$
F\,:\ [0,\infty)\ni I\mapsto \sigma I +r^{\frac{1}{q}}(x)I^{\frac{1-p}{q}}-\tilde{\kappa}
$$ is smooth, strictly increasing, with $F(0)=-\tilde\kappa\leq0$ and $\lim_{I\to\infty}F(I)=\infty$. Hence the existence and uniqueness of $I_{\sigma}(x,\tilde{\kappa})$  satisfying \eqref{LU2-1} follows from the intermediate value theorem. The continuity of $I_{\sigma}(x,\tilde{\kappa})$ with respect to $(x,\tilde{\kappa})\in \bar{\Om}\times[0,\infty)$ follows from the implicit function theorem.  It is easy to see from \eqref{LU2-1} that $I_{\sigma}(x,\tilde{\kappa})$ is strictly increasing in $\tilde{\kappa}$, $I_{\sigma}(x,0)\equiv 0$ and $\lim_{\tilde{\kappa}\to\infty}\min_{x\in\bar{\Om}}I_{\sigma}(x,\tilde{\kappa})=\infty$. Therefore, the function
$$
G\, : \ [0,\infty)\ni \tilde{\kappa}\mapsto \tilde{\kappa}|\Omega|+(1-\sigma)\int_{\Om}I_{\sigma}(s,\tilde{\kappa})dx
$$
is continuous and $G(0)=0$. Observe from \eqref{LU2-1} that
\begin{eqnarray*}
G(\tilde{\kappa})&=&\int_{\Om}(\tilde{\kappa}-\sigma I_{\sigma}(x,\tilde{\kappa}))dx+\int_{\Om}I_{\sigma}(x,\tilde{\kappa})dx\\
&=&\int_{\Om}r^{\frac{1}{q}}(x)[I_{\sigma}(x,\tilde{\kappa})]^{\frac{1-p}{q}}dx+\int_{\Om}I_{\sigma}(x,\tilde{\kappa})dx.
\end{eqnarray*}
Whence, $G$ is strictly increasing in $\tilde{\kappa}$ and $\lim_{\tilde{\kappa}\to\infty}G(\tilde{\kappa})=\infty$. We now invoke the intermediate value theorem to conclude that there is a unique $\tilde{\kappa}_{\sigma}>0$ satisfying $G(\tilde{\kappa}_{\kappa})=N$. This shows that $\tilde{\kappa}_{\sigma}$ is the unique positive  real number satisfying equation \eqref{LU2-2}.
\end{proof}

Based on Lemma \ref{lem-6}, we can prove Theorem \ref{TH2.4}(ii).

\begin{proof}[Proof of Theorem \ref{TH2.4}{\rm (ii)}] Recall that $\frac{\kappa}{d_S}$ satisfies \eqref{LU1} and \eqref{Eq6-7}. Hence, up to a subsequence, there is ${\kappa}^*\in [N_*,\frac{N}{|\Omega|}(1+\sigma)]$ such that
$$
\frac{\kappa}{d_S}\to {\kappa}^*\ \ \ \text{as}\  d_I+\big|\frac{d_I}{d_S}-\sigma\big|\to 0,
$$
where $N_*$ is given in \eqref{Eq6-7}. Note from \eqref{Eq6-1} that $I$ solves
\begin{equation}\label{LU2-3}
\begin{cases}
d_I\Delta I +\beta\Big[\Big(\frac{\kappa}{d_S}-\frac{d_I}{d_S}I\Big)^q-rI^{1-p}\Big]I^{p}=0,\ \ \ & x\in\Om,\cr
\p_{\nu}I=0, & x\in\p\Om.
\end{cases}
\end{equation}
By virtue of \eqref{Eq6-6} and the fact $S<\frac{\kappa}{d_S}$, we have
$$
\limsup_{d_I+\big|\frac{d_I}{d_S}-\sigma\big|\to0}\|I\|_{L^\infty(\Omega)}\le \Big[\Big(\frac{\beta}{\gamma}\Big)_{\max}\frac{N^q(1+\sigma)^q}{|\Omega|^q}\Big]^{\frac{1}{1-p}}.
$$
Hence, we can employ the singular perturbation theory for elliptic equations to assert from \eqref{LU2-3} that
$$
I\to I_{\sigma}(\cdot,{\kappa}^*)\ \ \ \ \ \text{uniformly on}\  \bar\Om\,\ \ \text{as}\ d_I+\big|\frac{d_I}{d_S}-\sigma\big|\to 0,
$$
where $I_{\sigma}(\cdot,{\kappa}^*)$ is given by Lemma \ref{lem-6}. Notice from Lemma \ref{Lem-2-1} that
$$
N=\int_{\Om}\Big[\Big(\frac{\kappa}{d_S}-\frac{d_I}{d_S}I\Big)+I\Big]
=\frac{\kappa}{d_S}|\Om|+\left(1-\frac{d_I}{d_S}\right)\int_{\Om}I.
$$
Taking $d_I+\big|\frac{d_I}{d_S}-\sigma|\to 0$, we  obtain
$$
N={\kappa}^*|\Om|+(1-\sigma)\int_{\Om}I_{\sigma}(x,{\kappa}^*).
$$
It then follows from Lemma \ref{lem-6} that ${\kappa}^*=\tilde{\kappa}_{\sigma}$. This shows that ${\kappa}^*$ is independent of the choice of subsequences. Therefore,
$$I\to I_{\sigma}(\cdot,\Tilde{\kappa}_{\sigma})\ \ \mbox{\ \ \,uniformly on $\bar\Om$\ \ as $d_I+\big|\frac{d_I}{d_S}-\sigma\big|\to 0$.}
$$
This together with \eqref{LU2-1} and $S=\frac{\kappa}{d_S}-\frac{d_I}{d_S}I$ implies that
$$
S\to r^{\frac{1}{q}}\Big[I_{\sigma}(\cdot,\tilde{\kappa}_{\sigma})\Big]^{\frac{1-p}{q}}\ \ \ \ \ \text{uniformly on}\ \bar\Om\,\ \ \text{as}\ d_I+\big|\frac{d_I}{d_S}-\sigma\big|\to 0.
$$
This establishes Theorem \ref{TH2.4}(ii).
\end{proof}

\section{Numerical simulations and conclusions}\label{Nuremical-Sim-section}
In this section, we first conduct numerical simulations to illustrate the results presented in this paper. Then we conclude the paper by providing a discussion on the biological implications of our findings.

\subsection{Simulations}
In this subsection, we specifically focus on the case of $p=1$ as we carry out numerical simulations to  illustrate the main results proved in this paper and complement certain cases where rigorous results have not been obtained yet. Biologically, $p=1$ means that the infection rate is proportional to the number of infected individuals. We set $q=0.5$. The domain $\Omega$ in our simulations is a unit circle centered at the origin.

\begin{figure}
     \centering
     \begin{subfigure}[b]{0.4\textwidth}
         \centering
          \includegraphics[scale=.3]{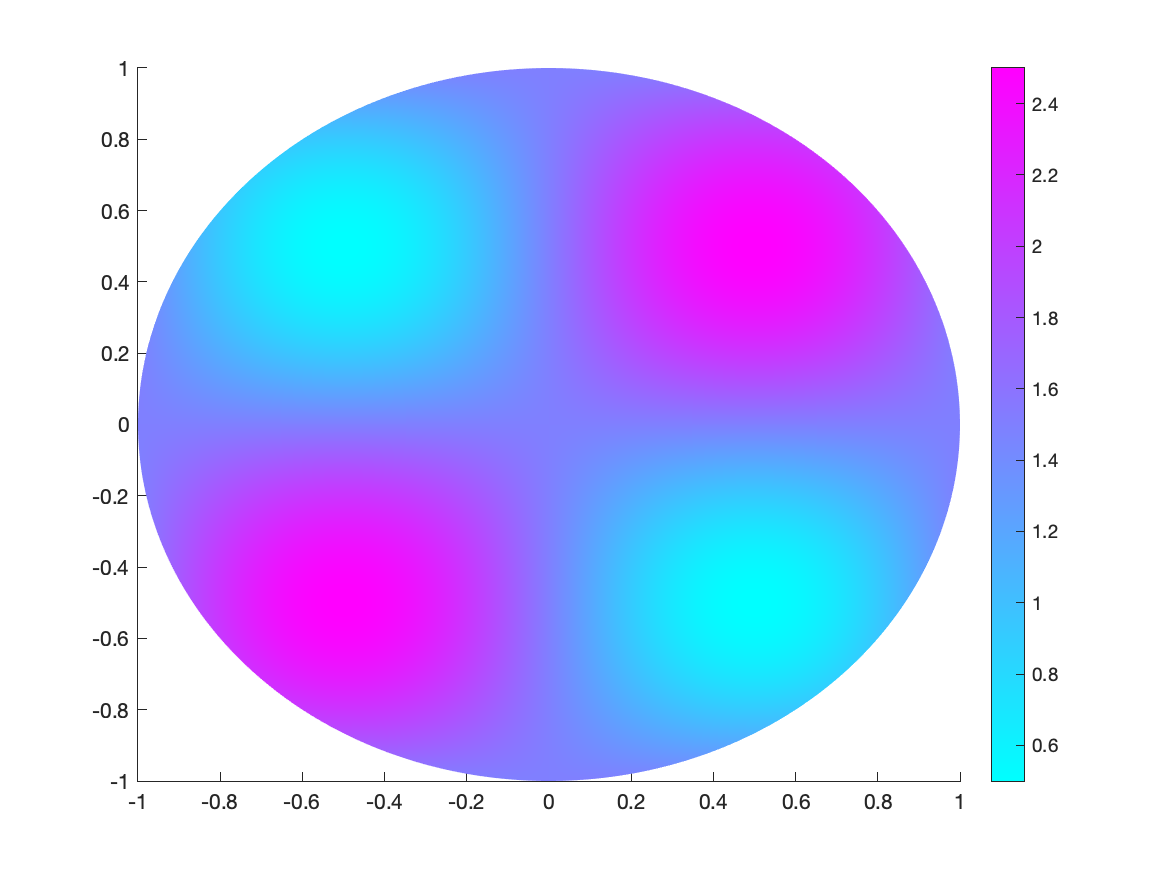}
         \caption{}
         \label{fig0:a1}
     \end{subfigure}
     \begin{subfigure}[b]{0.4\textwidth}
        \centering
     \includegraphics[scale=.3]{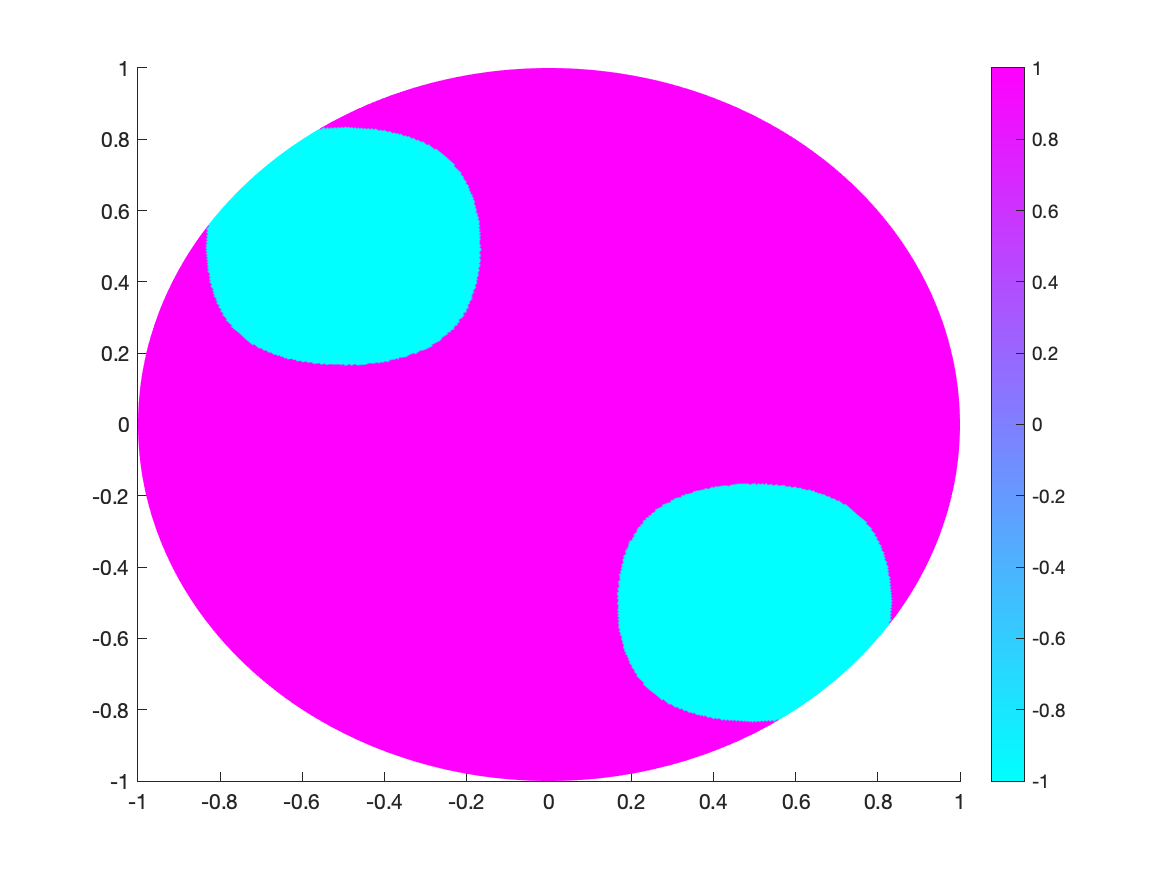}
             \caption{}
         \label{fig0:b1}
     \end{subfigure}
     \caption{Illustration of high risk ($R>1$) and low risk ($R<1$) domains when $p=1$ and $q=0.5$. (a): graph of $R(x, y)=N\beta(x, y)/(|\Omega|\gamma(x, y))$; (b): blue region is the low risk domain and red region is the high risk domain. The points $(0.5, 0.5)$ and $(-0.5, -0.5)$ are of the highest risk since the minimum of $r(x, y)=\gamma(x, y)/\beta(x, y)$ on $\bar\Omega$ is attained there.}
     \label{fig:0}
\end{figure}

\emph{Simulation 1}. The disease transmission and recovery rates are given by $\beta(x, y)=1.5+\sin(\pi x)\sin(\pi y)$ and $\gamma(x, y)=1$.  Figure \ref{fig:0} provides a visualization of the high and low risk domains. In Figure \ref{fig:0}(a), the function $R(x, y)=N\beta(x, y)/(|\Omega|\gamma(x, y))$ is depicted. If $R(x, y)<1$, the point $(x, y)$ is considered to be of lower risk; if $R(x, y)>1$, it is deemed to be of higher risk. The lower risk domains are situated in the upper left and bottom right regions of $\Omega$, as shown in Figure \ref{fig:0}(b). Additionally, the minimum of $r(x, y)=\gamma(x, y)/\beta(x, y)$ on $\bar\Omega$ is attained at the points $(0.5, 0.5)$ and $(-0.5, -0.5)$, which are the points of the highest risk. Since the equilibrium problem always has a trivial solution, we will solve the reaction-diffusion system \eqref{model-1} instead of the elliptic system \eqref{1-1}.
For our simulation, we initialize the system with $S_0=0.8$ and $I_0=0.2$. Then it holds that $\frac{N}{|\Omega|}>r^\frac{1}{q}_{\min}$. By Lemma \ref{lem_R0} and Theorem \ref{theorem_ex}, model \eqref{1-1} has at least one EE solution for any $d_S, d_I>0$.

\begin{figure}
     \centering
     \begin{subfigure}[b]{0.48\textwidth}
         \centering
          \includegraphics[scale=.3]{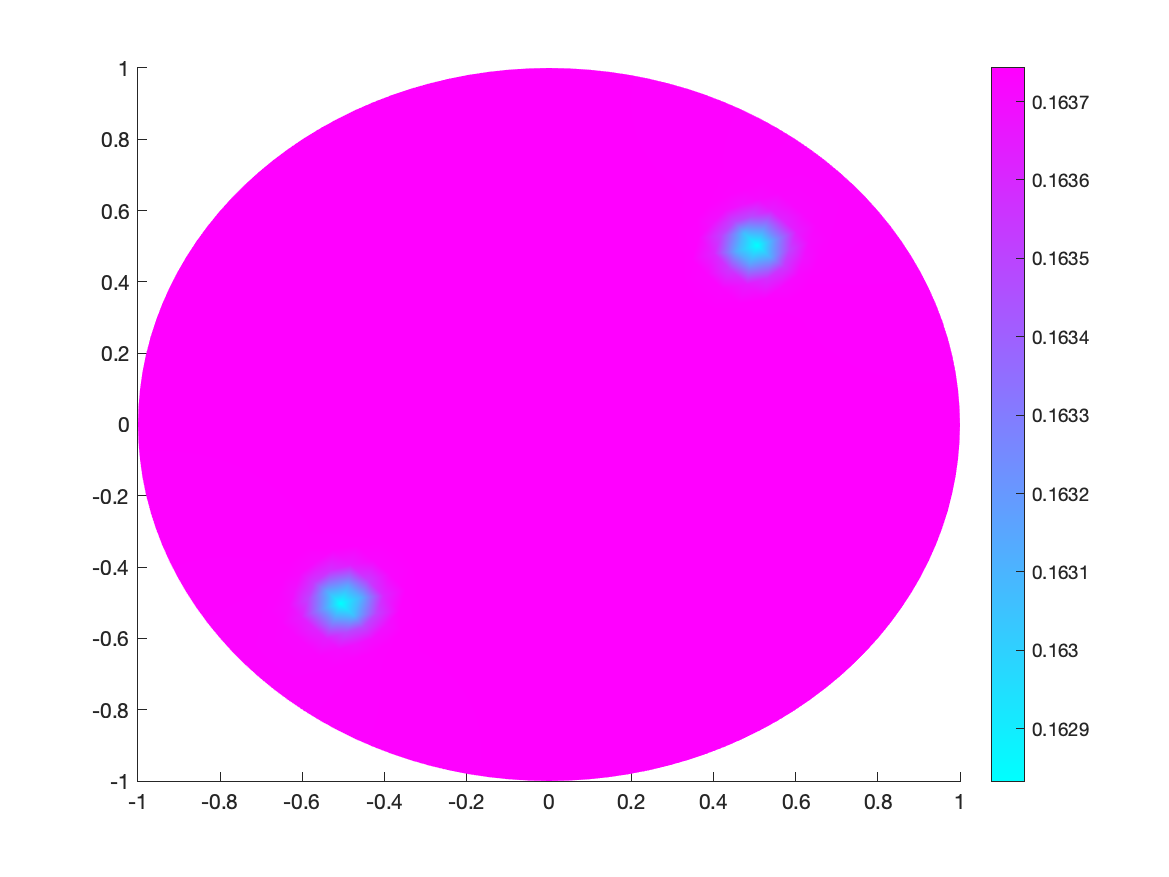}
         \caption{$S(x, y, 200)$ when $d_S=1, d_I=10^{-5}$.}
         \label{fig:a1}
     \end{subfigure}
     \begin{subfigure}[b]{0.48\textwidth}
        \centering
     \includegraphics[scale=.3]{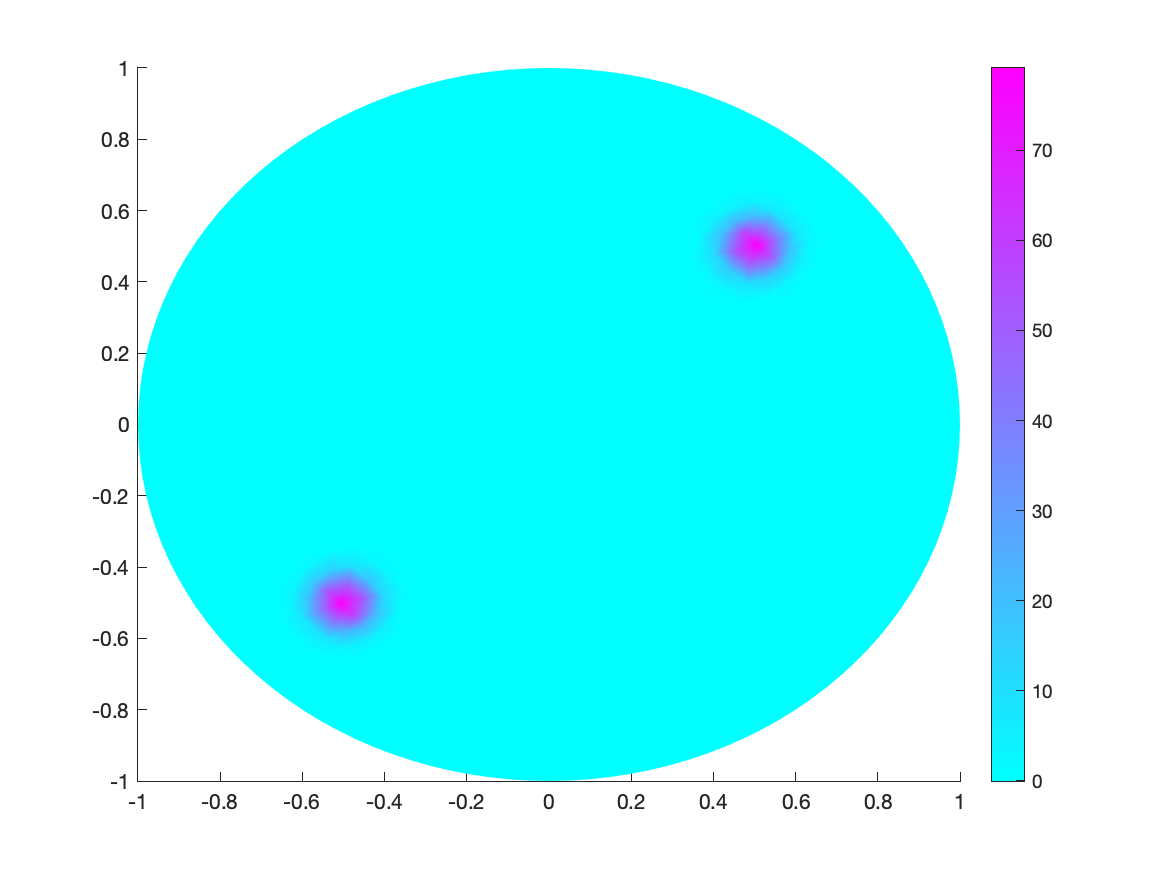}
             \caption{$I(x, y, 200)$ when $d_S=1, d_I=10^{-5}$.}
         \label{fig:b1}
     \end{subfigure}

\vskip20pt
         \begin{subfigure}[b]{0.48\textwidth}
         \centering
          \includegraphics[scale=.3]{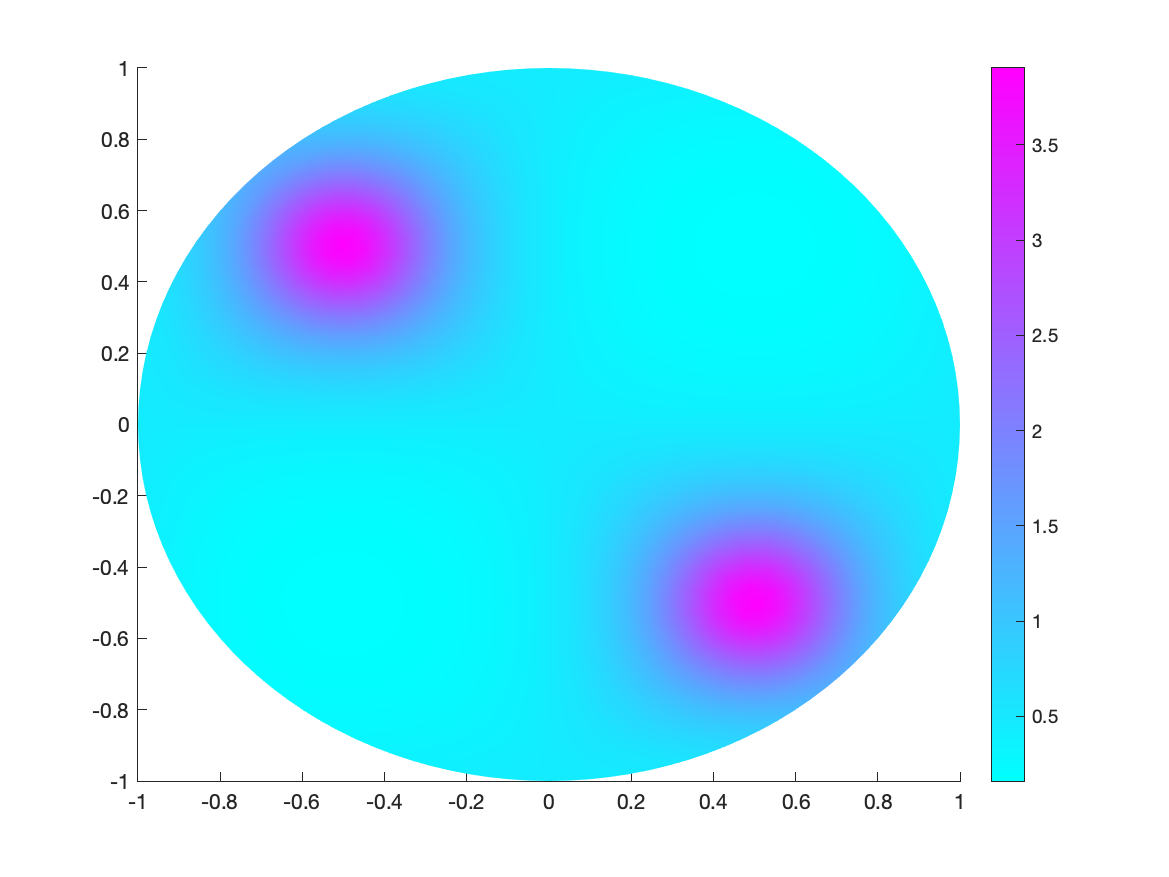}
         \caption{$S(x, y, 200)$ when $d_S=10^{-5}, d_I=1$.}
         \label{fig:a2}
     \end{subfigure}
     \begin{subfigure}[b]{0.48\textwidth}
        \centering
     \includegraphics[scale=.3]{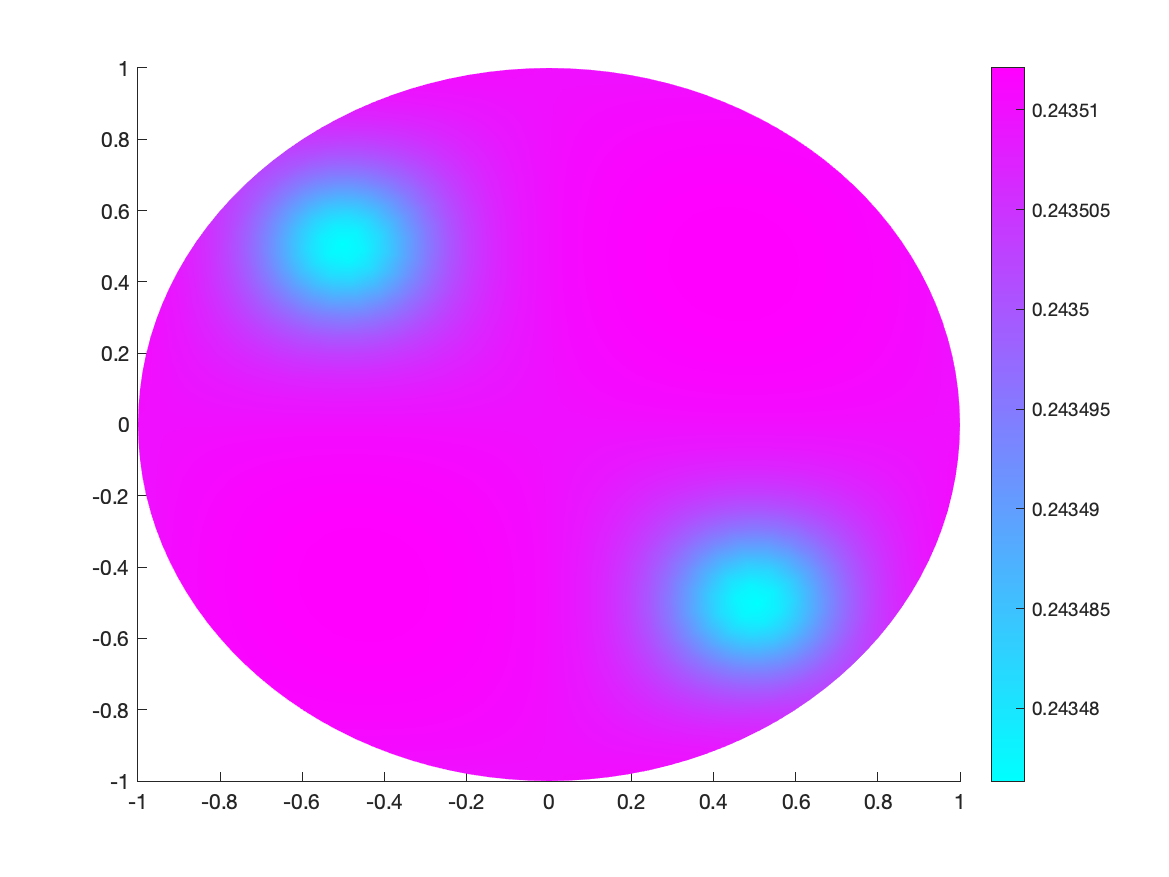}
             \caption{$I(x, y, 200)$ when $d_S=10^{-5}, d_I=1$.}
         \label{fig:b2}
     \end{subfigure}

\vskip20pt
         \begin{subfigure}[b]{0.48\textwidth}
         \centering
          \includegraphics[scale=.3]{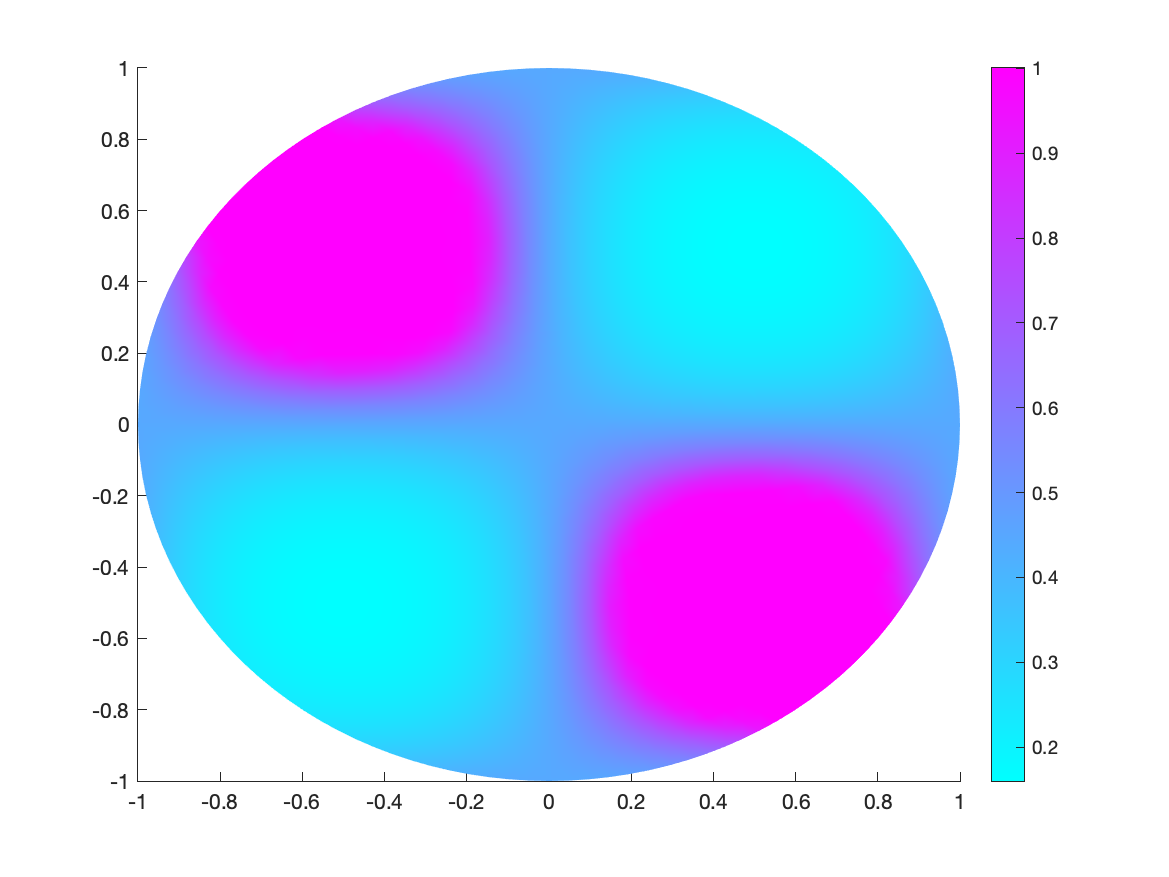}
         \caption{$S(x, y, 200)$ when $d_S= d_I=10^{-5}$.}
         \label{fig:a3}
     \end{subfigure}
     \begin{subfigure}[b]{0.48\textwidth}
        \centering
     \includegraphics[scale=.3]{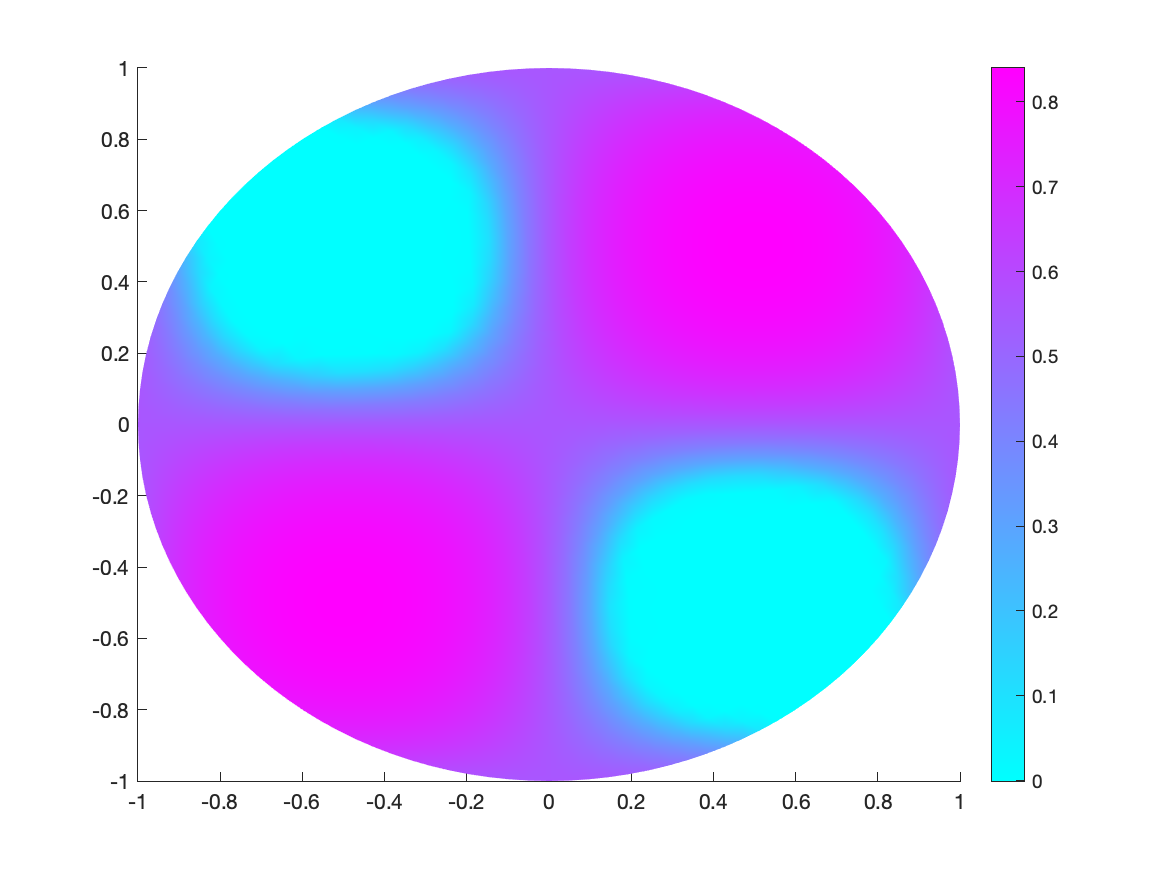}
             \caption{$I(x, y, 200)$ when $d_S= d_I=10^{-5}$.}
         \label{fig:b3}
     \end{subfigure}
     \caption{Simulation of \eqref{model-1} with $p=1$ and $q=0.5$. The disease transmission and recovery rates are $
\beta(x, y)=1.5+\sin(\pi x)\sin(\pi y)$ and $\gamma(x, y)=1$, respectively.}
     \label{fig:1}
\end{figure}

When we restrict the movement of infected individuals only ($d_S=1, d_I=10^{-5}$), as depicted in Figures \ref{fig:1}(a)-(b), the infected individuals tend to concentrate at $(0.5, 0.5)$ and $(-0.5, -0.5)$, which are precisely the points of the highest risk. Additionally, the distribution of susceptible individuals is nearly spatially homogeneous. This aligns with the findings of Theorem \ref{TH2.1}(i).

On the other hand, if we limit the movement of susceptible individuals only ($d_S=10^{-5}, d_I=1$), as shown in Figures \ref{fig:1}(c)-(d), the infected individuals tend to be spatially evenly distributed, and there are more susceptible individuals residing in the lower risk region compared to the higher risk region. This observation corresponds to the predictions of Theorem \ref{TH2.3}(i-c).

When we restrict the movement of both susceptible and infected individuals ($d_S= d_I=10^{-5}$), as illustrated in Figures \ref{fig:1}(e)-(f), it appears that the infection in some lower risk regions is eliminated, and the density of susceptible individuals is higher in lower risk regions than in higher risk regions. These observations validate the claims made in Theorem \ref{TH2.4}(i).

\begin{figure}
     \centering
         \begin{subfigure}[b]{0.49\textwidth}
         \centering
          \includegraphics[scale=.3]{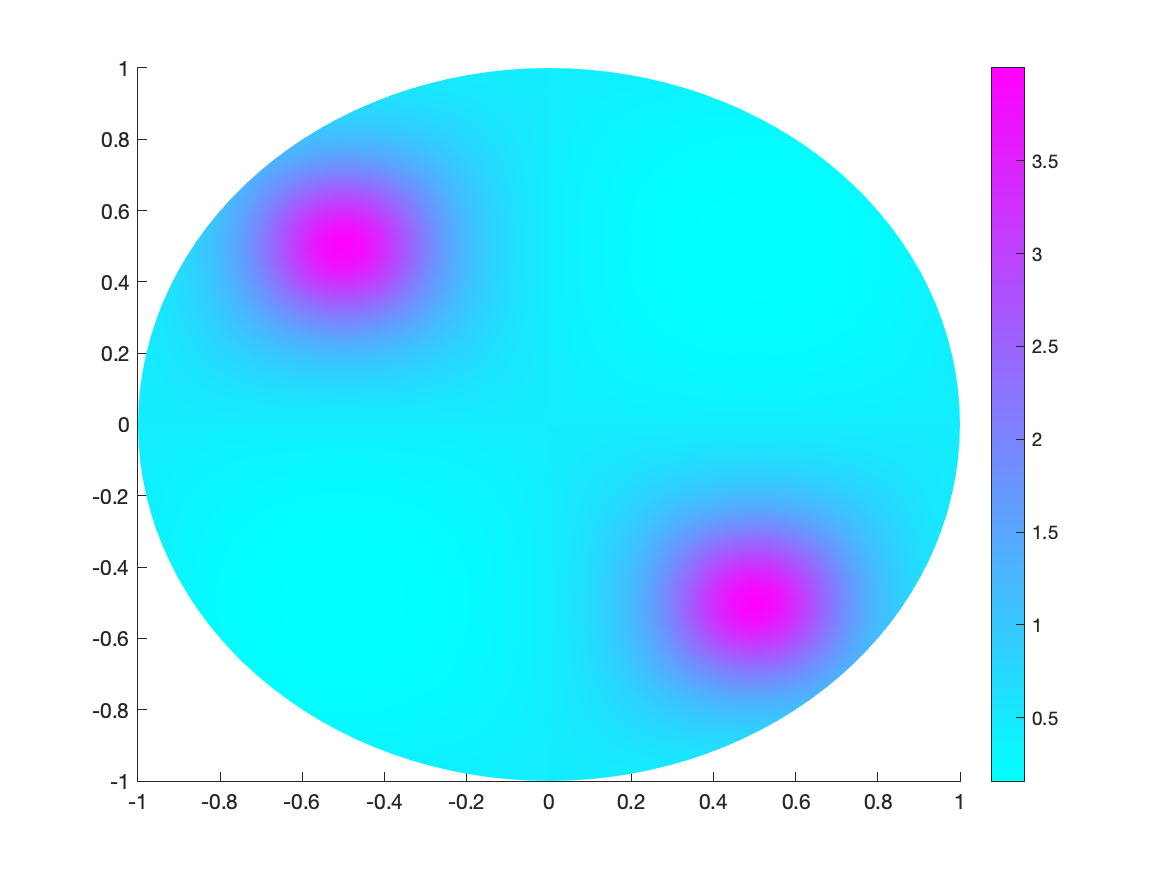}
         \caption{$S(x, y, 2000)$ when $d_S=10^{-7}, d_I=10^{-3}$.}
         \label{fig2a}
     \end{subfigure}
     \begin{subfigure}[b]{0.48\textwidth}
        \centering
     \includegraphics[scale=.3]{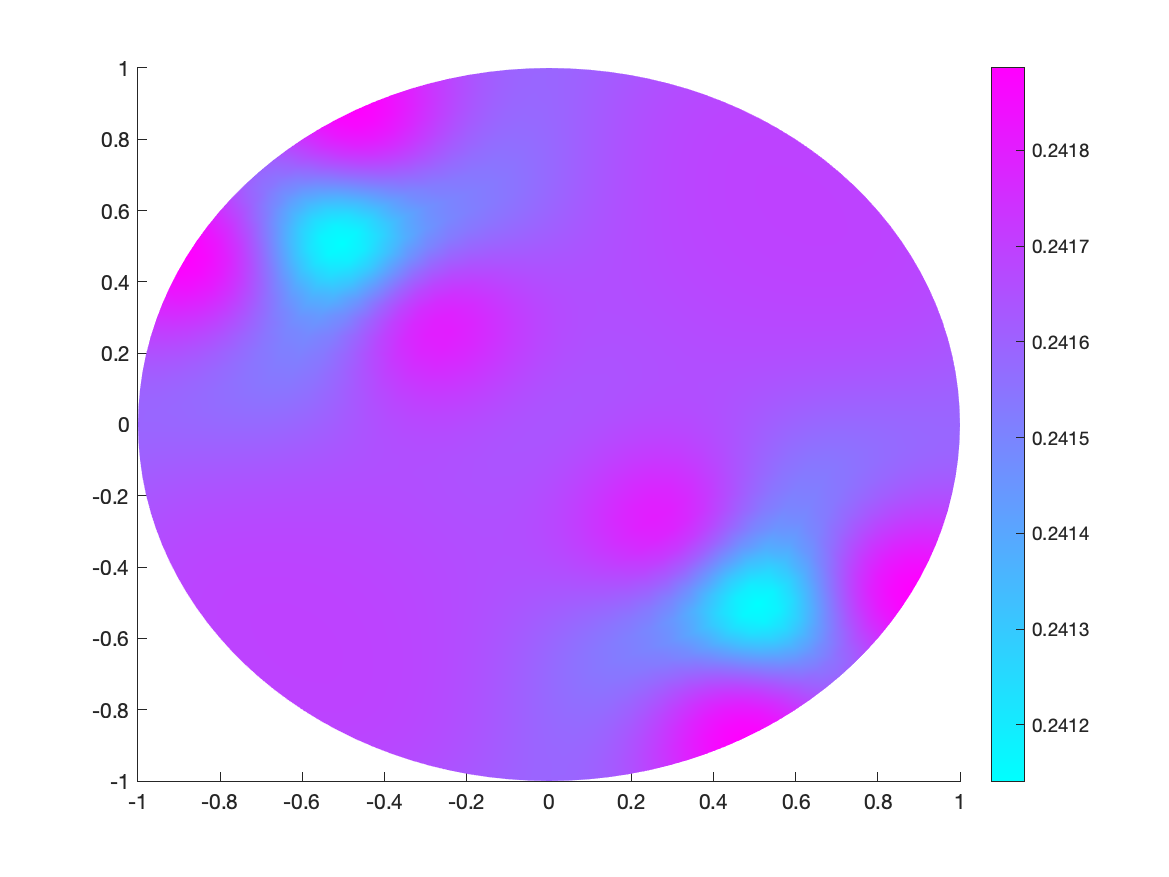}
             \caption{$I(x, y, 2000)$ when $d_S=10^{-7}, d_I=10^{-3}$.}
         \label{fig2b}
     \end{subfigure}

\vskip25pt

   \begin{subfigure}[b]{0.48\textwidth}
         \centering
          \includegraphics[scale=.3]{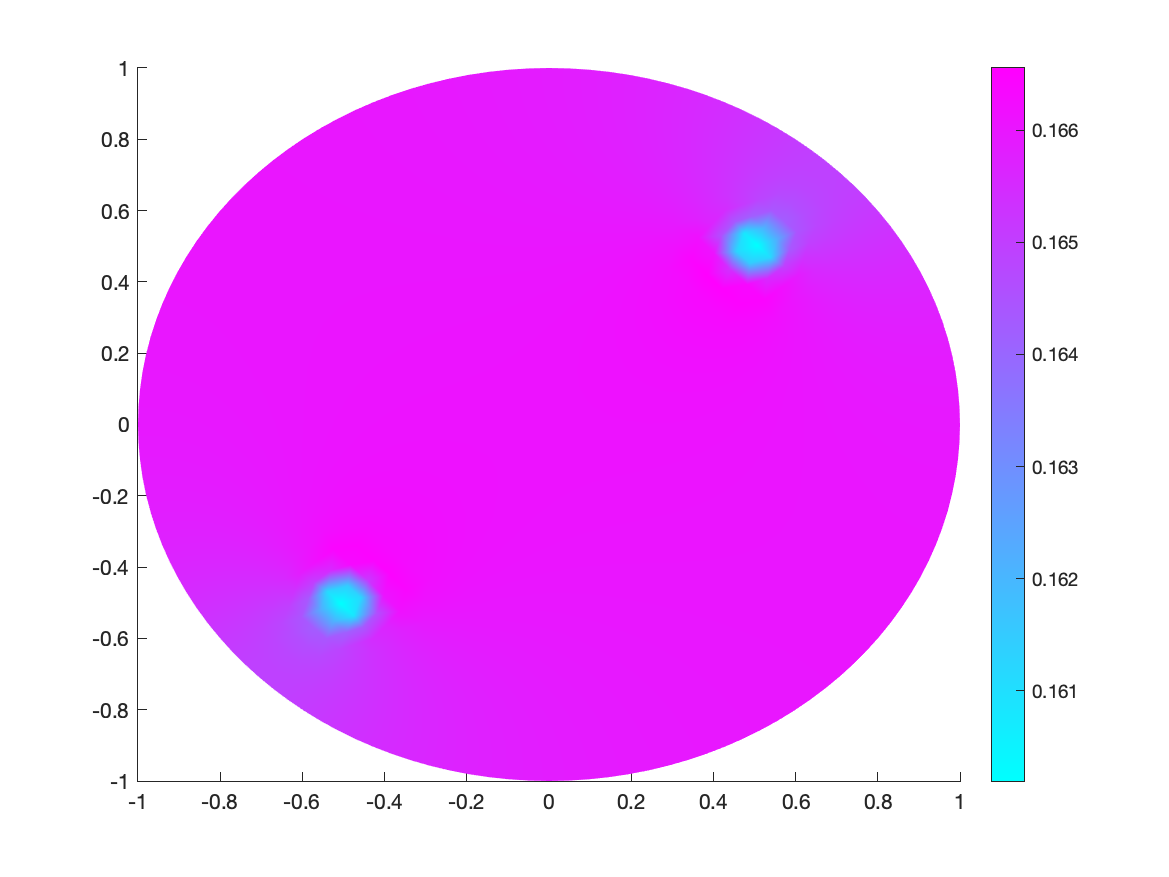}
         \caption{$S(x, y, 2000)$ when $d_S=10^{-2}, d_I=10^{-7}$.}
         \label{fig2c}
     \end{subfigure}
     \begin{subfigure}[b]{0.48\textwidth}
        \centering
     \includegraphics[scale=.3]{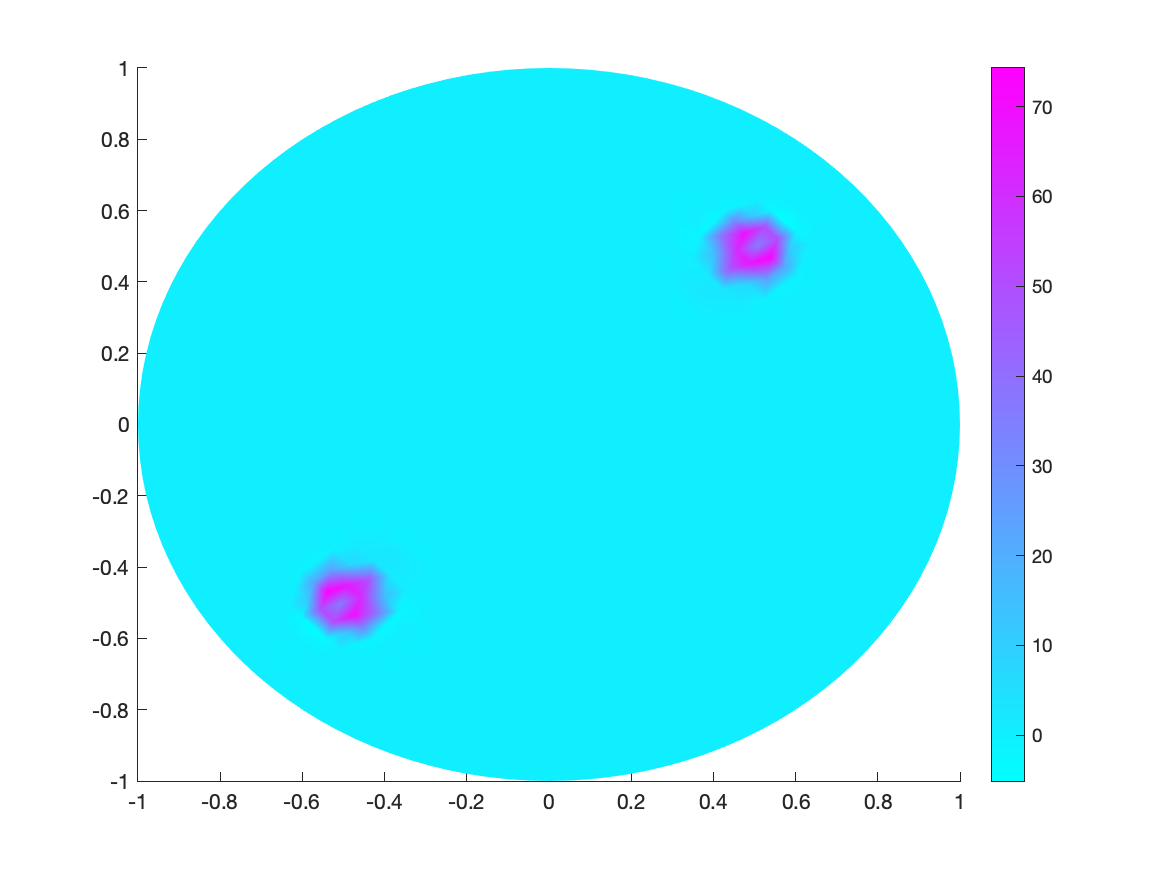}
             \caption{$I(x, y, 2000)$ when $d_S=10^{-2}, d_I=10^{-7}$.}
         \label{fig2d}
     \end{subfigure}

\vskip25pt
         \begin{subfigure}[b]{0.48\textwidth}
         \centering
          \includegraphics[scale=.3]{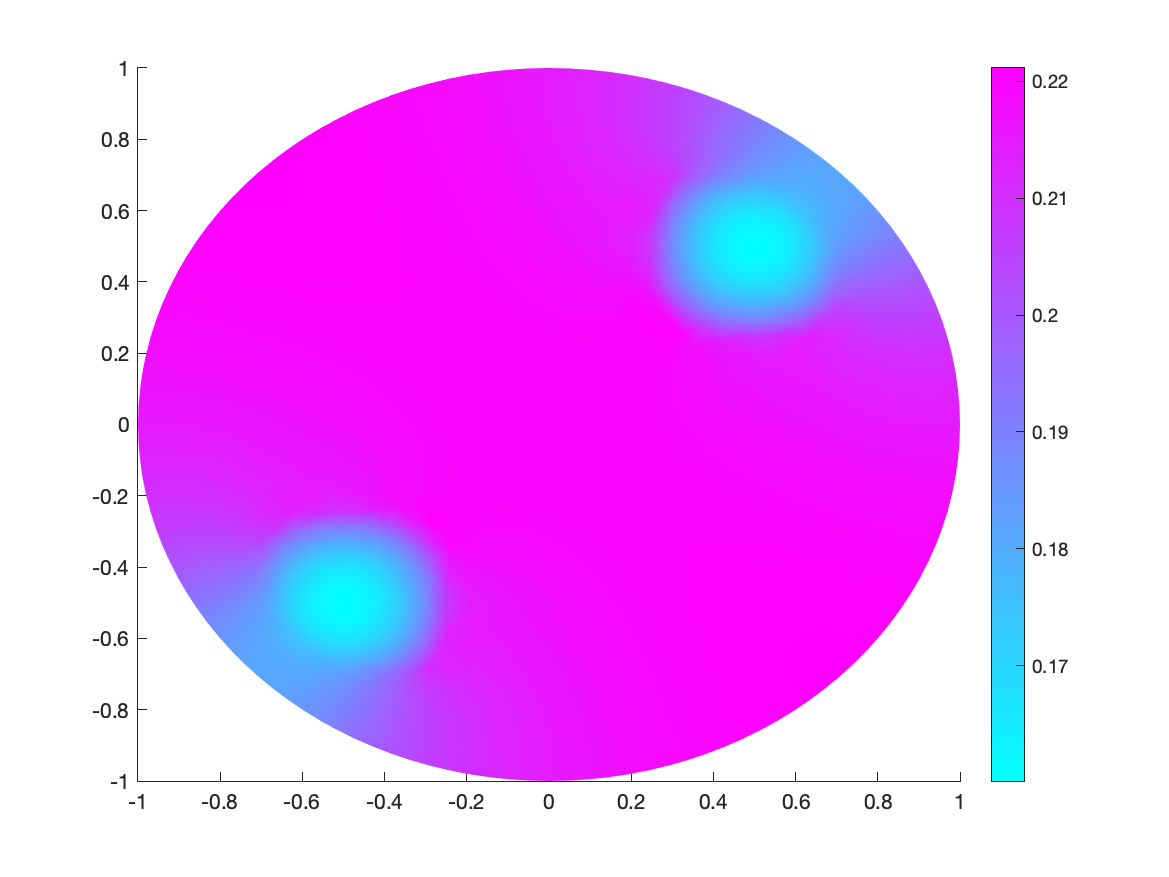}
         \caption{$S(x, y, 2000)$ when $d_S=10^{-3}, d_I=10^{-7}$.}
         \label{fig2e}
     \end{subfigure}
     \begin{subfigure}[b]{0.48\textwidth}
        \centering
     \includegraphics[scale=.3]{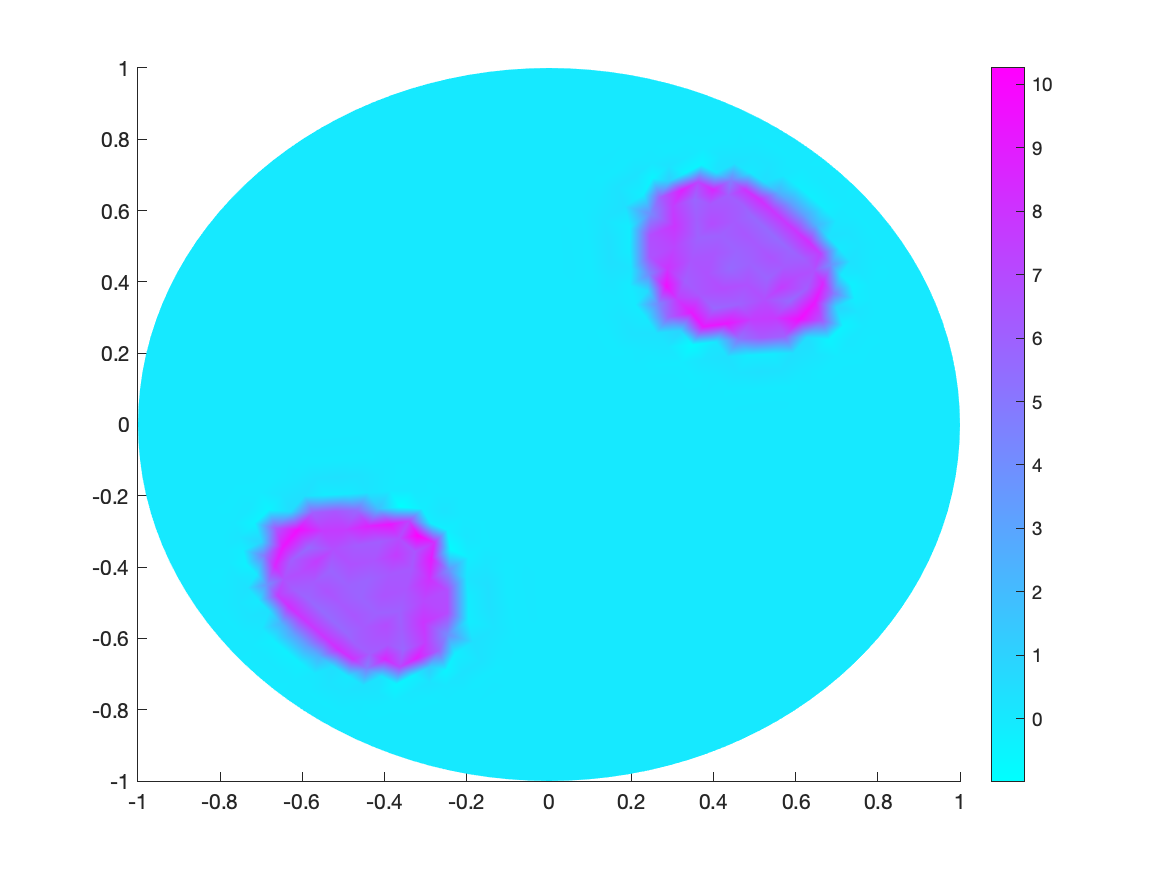}
             \caption{$I(x, y, 2000)$ when $d_S=10^{-3}, d_I=10^{-7}$.}
         \label{fig2f}
     \end{subfigure}

     \caption{Simulation of \eqref{model-1} with $p=1$ and $q=0.5$. The disease transmission and recovery rates are $
\beta(x, y)=1.5+\sin(\pi x)\sin(\pi y)$ and $\gamma(x, y)=1$, respectively.}
     \label{fig:2}
\end{figure}

Finally, we restrict the movement of  susceptible and infected individuals disproportionately. If the movement rate of the susceptible is considerably smaller than that of the infected ($d_S=10^{-7},  d_I=10^{-3}$), Figure \ref{fig2b} shows that the density of the infected people converges to a constant ($\approx$ 0.24). This aligns with Theorem \ref{TH2.4}(i-c), which predicts that the EEs converge to $\left(r^{\frac{1}{q}},\ \frac{1}{|\Omega|}\Big(N-\int_{\Omega}r^{\frac{1}{q}}\Big)\right)$ in this case. Note that $\ \frac{1}{|\Omega|}\Big(N-\int_{\Omega}r^{\frac{1}{q}}\Big)\approx 0.24$.

If the movement rate of the infected is considerably smaller than that of the susceptible ($d_S=10^{-2},  d_I=10^{-7}$), Figures \ref{fig2c}-\ref{fig2d} show that the density of the susceptible converges to a  constant ($\approx$ 0.16) and the infected people concentrate near the two points of the highest risk. This agrees with Theorem \ref{TH2.4}(i-d).

We remark that when the movement rates of both susceptible and infected individuals are small,  they may need to fall into some ranges for us to observe the phenomenon in Figures \ref{fig2a}-\ref{fig2d}. For example, if we set $d_S=10^{-3}$ and  $d_I=10^{-7}$, Figures  \ref{fig2e}-\ref{fig2f} show that the infected people are aggregating around the high risk domain at $t=2000$, and we are not clear whether the concentration phenomenon in Figures \ref{fig2c}-\ref{fig2d} appears in this case.

\emph{Simulation 2}. The disease transmission and recovery rates are
$$
\beta(x, y)=0.5\ \ \text{and} \ \gamma(x, y)=f(x)f(y),
$$
respectively, where
\begin{equation}\label{f}
    f(x)=
    \left\{
\begin{array}{ll}
  0.5+0.4x^2,   & x\le 0, \\
  0.5,   & 0<x\le 0.25,\\
  0.5+0.4(x-0.25)^2, & 0.25<x\le 0.5,\\
  0.5+1.6(x-0.625)^2,& 0.5<x\le 1.
\end{array}
\right.
\end{equation}
 The minimum of $r(x, y)=\gamma(x, y)/\beta(x, y)$ on $\bar\Omega$ is attained at
 $$
 \Omega^*=[0, 0.25]\times [0, 0.25]\cup (\{0.625\}\times [0, 0.25])\cup  ([0, 0.25]\times \{0.625\})\cup \{(0.625, 0.625)\}.
 $$
 Then, $\Omega^*$ is the set of points of the highest risk, and it consists with a rectangular region, two line segments, and one isolated point.

\begin{figure}
     \centering
     \begin{subfigure}[b]{0.48\textwidth}
         \centering
          \includegraphics[scale=.3]{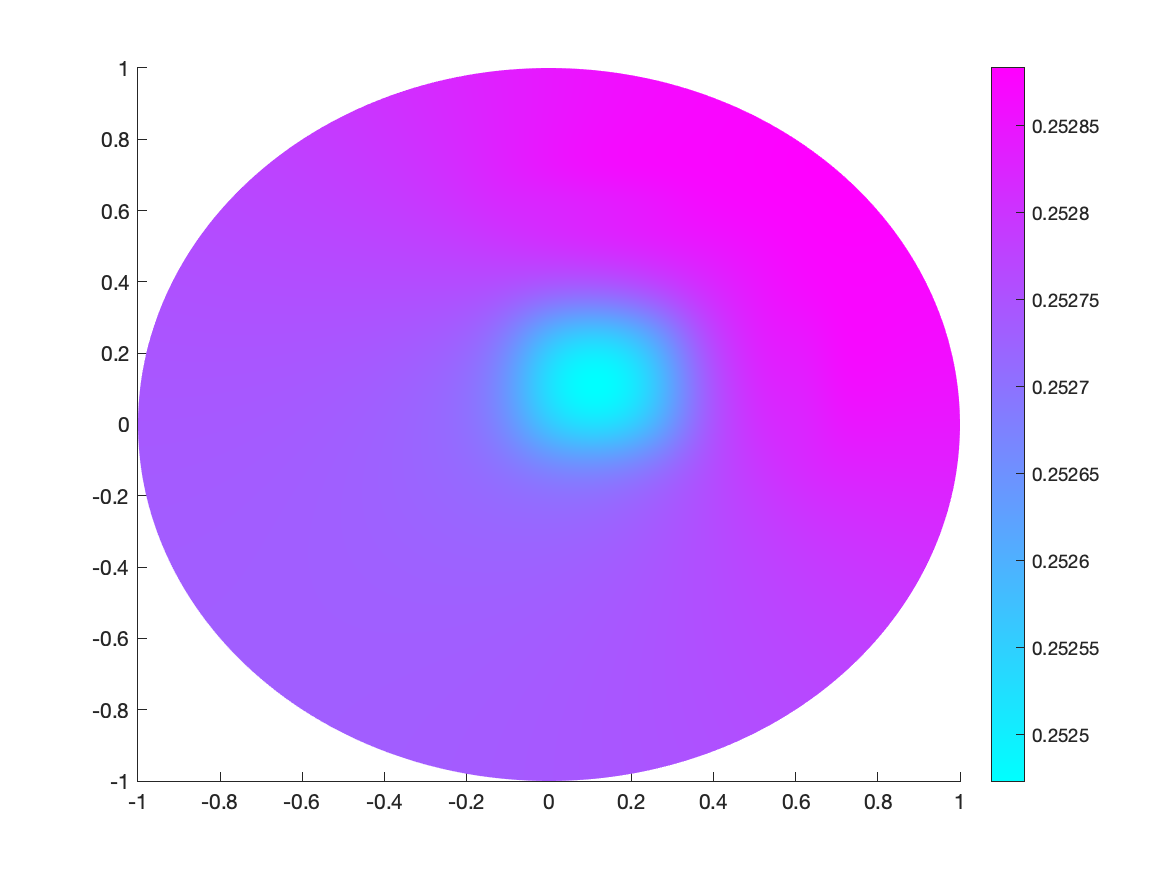}
        \caption{$S(x, y, 800)$ when $d_S=1$ and $d_I=10^{-5}$.}
         \label{fig:a14}
     \end{subfigure}
     \begin{subfigure}[b]{0.48\textwidth}
        \centering
     \includegraphics[scale=.3]{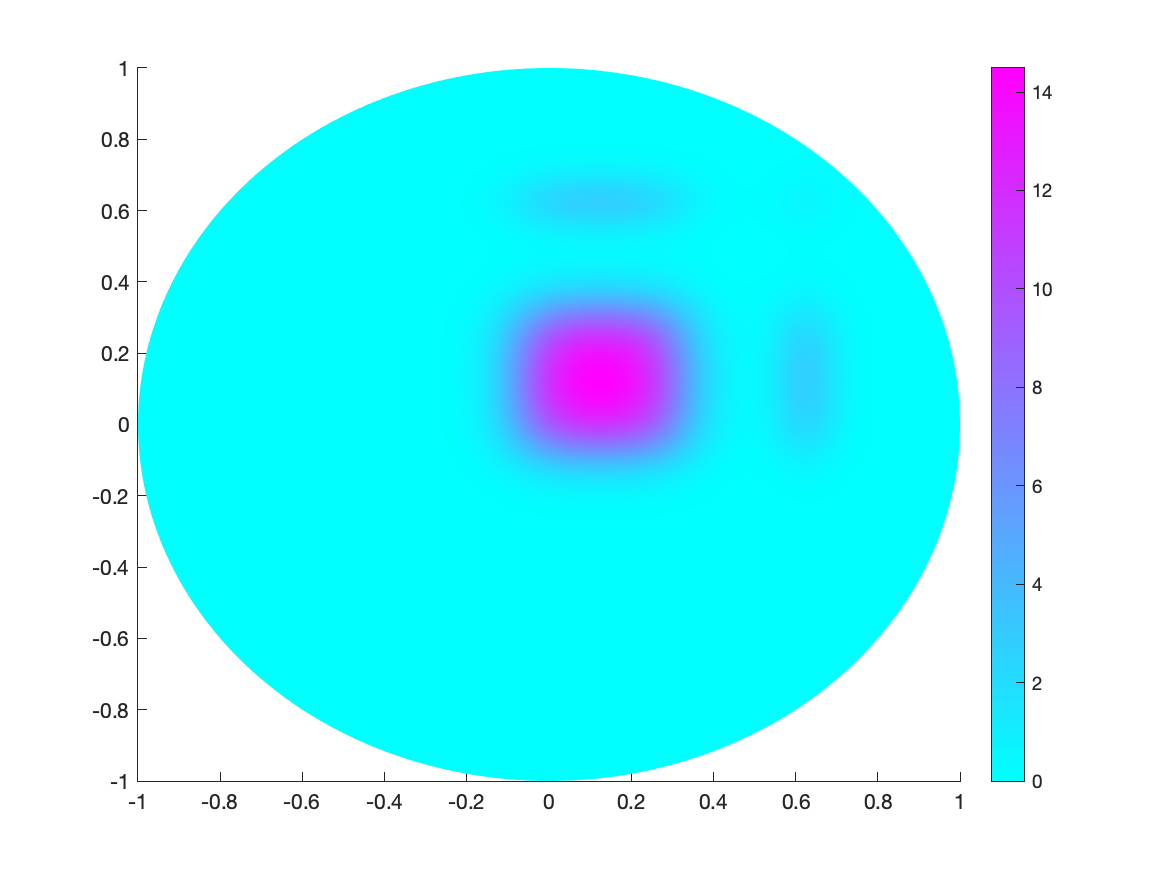}
             \caption{$I(x, y, 800)$ when $d_S=1$ and $d_I=10^{-5}$.}
         \label{fig:b14}
     \end{subfigure}

               \begin{subfigure}[b]{0.48\textwidth}
         \centering
          \includegraphics[scale=.3]{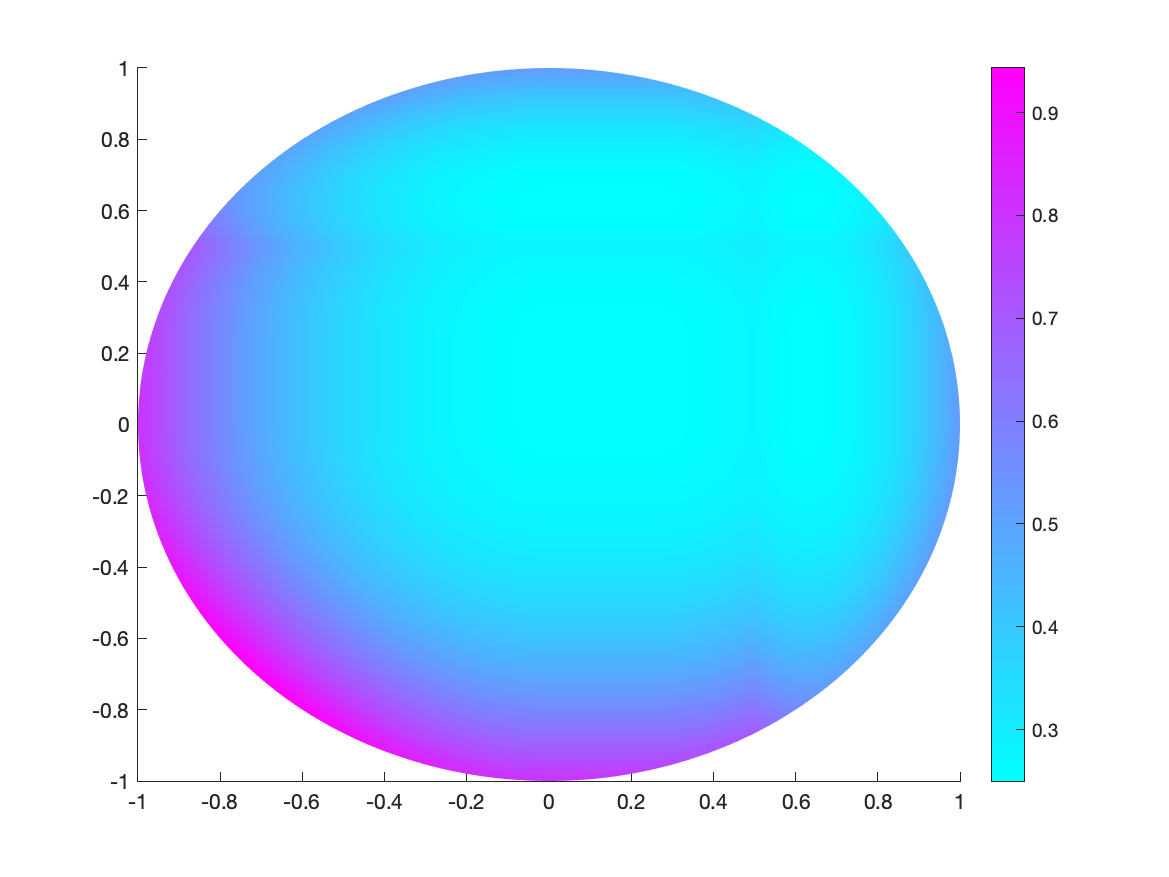}
        \caption{$S(x, y, 800)$ when $d_S=10^{-5}$ and $d_I=1$.}
     \end{subfigure}
     \begin{subfigure}[b]{0.48\textwidth}
        \centering
     \includegraphics[scale=.3]{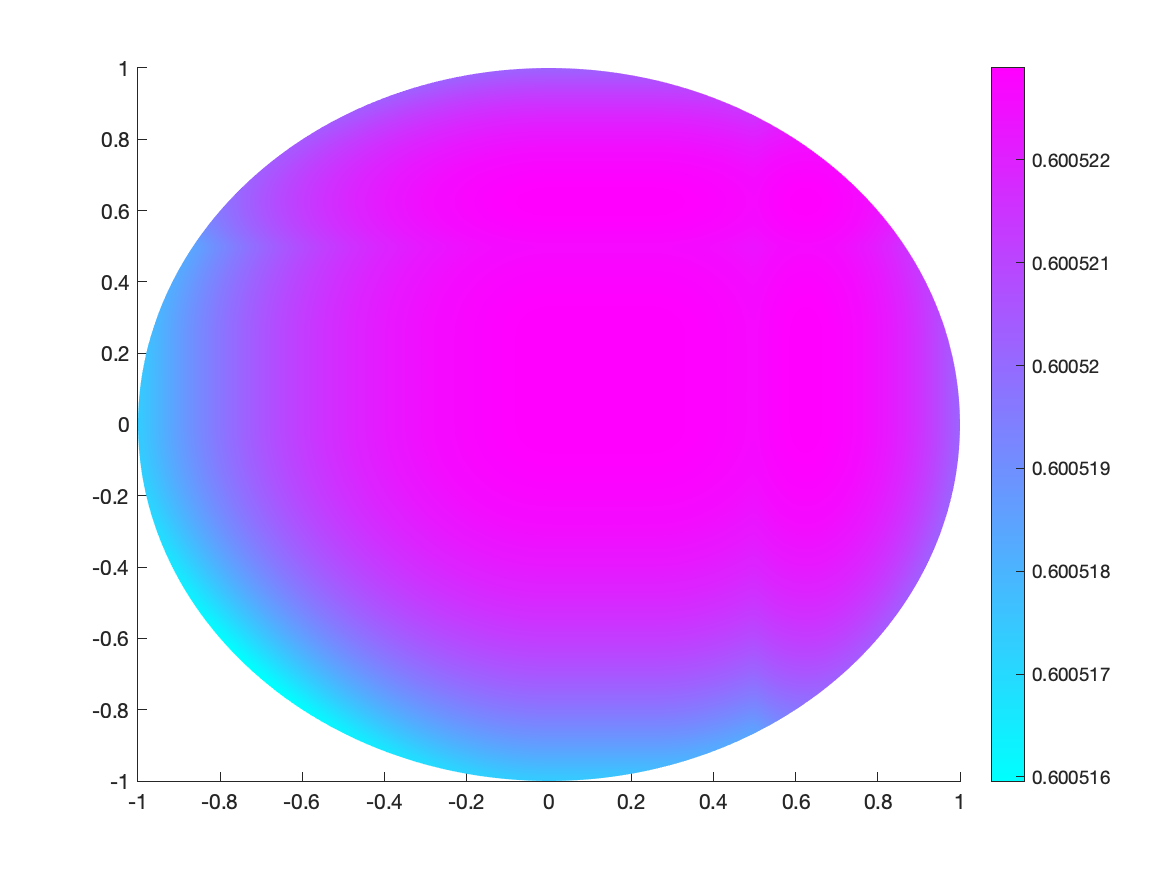}
             \caption{$I(x, y, 800)$ when $d_S=10^{-5}$ and $d_I=1$.}
         \label{fig:b15}
     \end{subfigure}

    \begin{subfigure}[b]{0.49\textwidth}
         \centering
          \includegraphics[scale=.3]{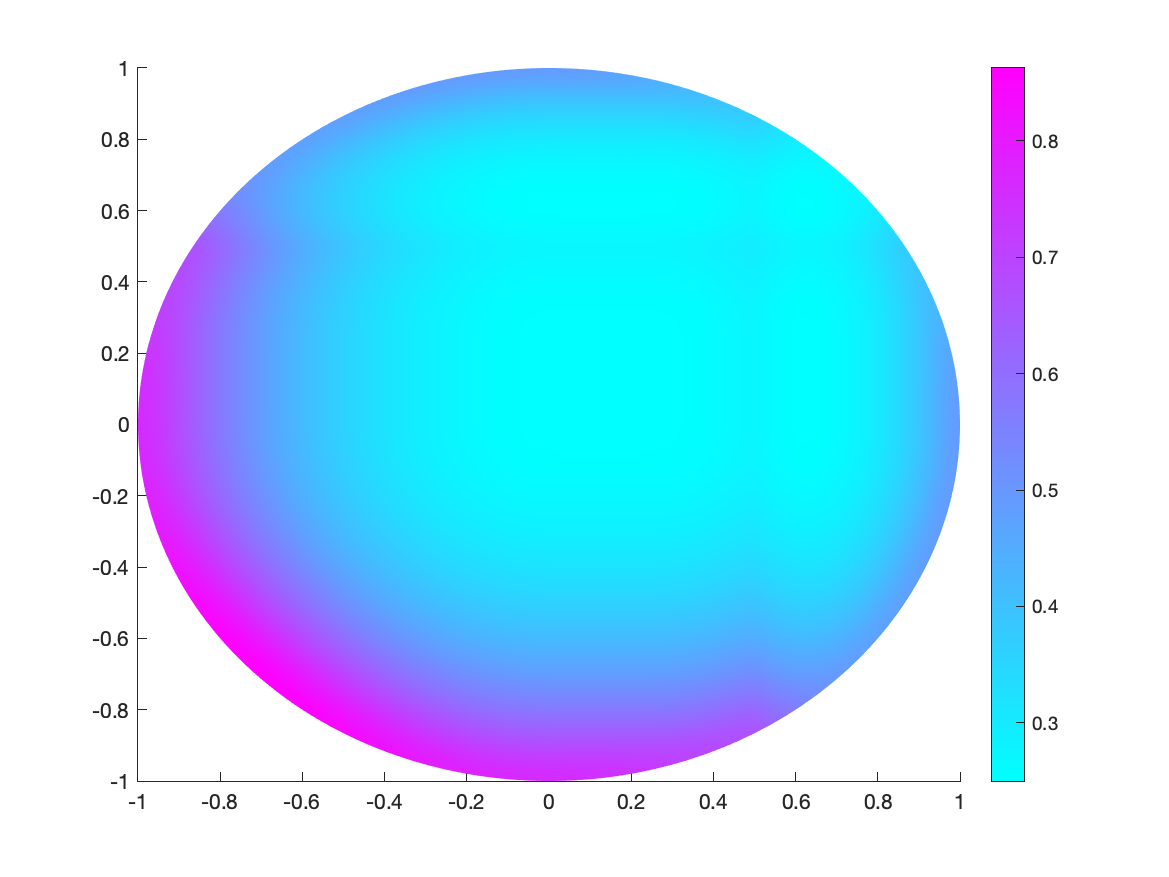}
        \caption{$S(x, y, 2000)$ when $d_S=d_I=10^{-5}$.}
     \end{subfigure}
     \begin{subfigure}[b]{0.49\textwidth}
        \centering
     \includegraphics[scale=.3]{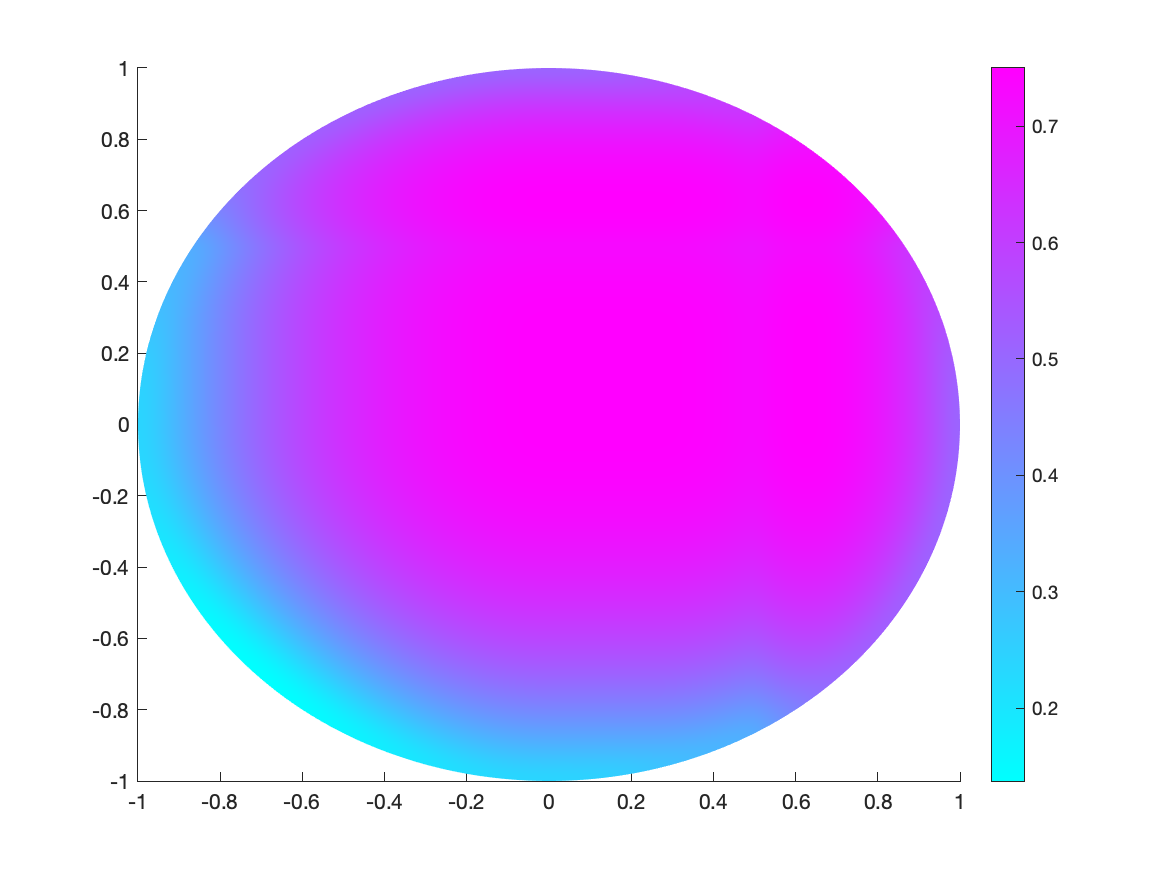}
             \caption{$I(x, y, 2000)$ when $d_S=d_I=10^{-5}$.}
         \label{fig:b16}
     \end{subfigure}

     \caption{Simulation of model \eqref{model-1} with $p=1$, and $q=0.5$. The disease transmission and recovery rates are $\beta(x, y)=0.5$ and $\gamma(x, y)=f(x)f(y)$ respectively, where $f$ is given by \eqref{f}. The minimum of $r=\gamma/\beta$ is attained at  $[0, 0.25]\times [0, 0.25]\cup (\{0.625\}\times [0, 0.25])\cup  ([0, 0.25]\times \{0.625\})\cup \{(0.625, 0.625)\}$, which consist with the points of the highest risk.}
     \label{fig:5}
\end{figure}

 When restricting the movement of infected people  ($d_S=1, d_I=10^{-5}$), Figure \ref{fig:5}(b) shows the infected people will reside in the rectangular region $[0, 0.25]\times [0, 0.25]$ and line segments $\{0.625\}\times [0, 0.25])\cup  ([0, 0.25]\times \{0.625\}$. However, there seems to be no infected people around the isolated highest risk point $\{(0.625, 0.625)\}$, which has not been reflected in our theoretical results. Indeed, in one-dimensional domain, this has been confirmed in \cite{PWZZ2023} both rigorously and numerically. Figure \ref{fig:5}(a) shows that the susceptible people distribute evenly in the domain. These figures are consistent with Theorem \ref{TH2.1}.

When restricting the movement of susceptible people ($d_S=10^{-5}, d_I=1$), Figure \ref{fig:5}(d) shows that the infected people distribute evenly in the domain, which is consistent with Theorem \ref{TH2.3}. When restricting the movement of both susceptible and infected people ($d_S= d_I=10^{-5}$), as shown in Figure \ref{fig:5}(f), the infection cannot be eliminated in any region of the domain. This is different from Figure \ref{fig:1}(f), where the infection in some part of the region is eliminated. This is still consistent with  Theorem \ref{TH2.4}(i), which claims that $(S, I)\to (S_1, I_{1}):=(\min\{\tilde{\kappa}_{1},r^{\frac{1}{q}}\}, \Big(\tilde{\kappa}_{1}-r^{\frac{1}{q}}\Big)_+)$ as $d_I+|\frac{d_I}{d_S}-1|\to 0$. In this case, we have $\tilde\kappa_1>r^\frac{1}{q}$ and $(S, I)\to (r^\frac{1}{q}, \tilde{\kappa}_{1}-r^{\frac{1}{q}})$. Indeed, Figure \ref{fig:5}(e) shows that $S\to r^\frac{1}{q}$.

\subsection{Conclusions}

In this subsection, we discuss our results on the asymptotic profiles of the EEs of model \eqref{model-1} from the perspective of disease control. We first suppose that the infection rate is proportional to the number of infected individuals, i.e., $p=1$.

The results for the asymptotic profiles of the EE as $d_I\to 0$ (Theorem \ref{TH2.1} and Proposition \ref{Prop2.1}) focus on the strategy of limiting the movement of infected individuals.
  Theorem \ref{TH2.1}(i) indicates that the infected individuals will concentrate on some, if not all, of the highest-risk points, if the highest-risk habitat $\Omega_*$ consists of finitely many isolated points. This result is confirmed by our numerical simulations (see Figures \ref{fig:1}(a)-(b)). Proposition \ref{Prop2.1} suggests that if $\Omega_*$ contains at least one non-empty subdomain, no concentration phenomenon can occur, and the disease will aggregate in some, if not all, of these highest-risk subdomains. However, if the highest-risk habitat $\Omega_*$ satisfies $|\Omega_*|=0$ and contains a submanifold, it remains an open problem whether the disease can also concentrate on such a submanifold.

The results for the asymptotic profiles of the EE as $d_S\to 0$ (Theorems \ref{TH2.3}) focus on the strategy of limiting the movement of susceptible individuals.  Theorem \ref{TH2.3}(i) indicates that the outcome of the control strategy depends on the population size: if the total population size is small, i.e., $N\le \int_{\Omega}r^\frac{1}{q}$, the disease will be eliminated; if the size is large, i.e., $N>\int_{\Omega}r^\frac{1}{q}$, the disease cannot be controlled. Figures \ref{fig:1}(c)-(d) presents a situation where the total population size is large, and the disease is not controlled.

The results for the asymptotic profiles of the EE as $d_I+|\frac{d_I}{d_S}-\sigma|\to 0$ (Theorems \ref{TH2.4}) are for the strategy of limiting the movement of both susceptible and infected individuals at some proportional rate. Theorem \ref{TH2.4} (i) indicates that the infection will be eliminated only in some region of the domain if the infection rates are proportionately restricted, which is supported by Figures \ref{fig:1}(e)-(f). If the movement rate of the susceptible people is much smaller (i.e. $\sigma\to\infty$) and the total population is large, the infection cannot be eliminated and the density of the infected people converges to a constant (see Figures \ref{fig2a}-\ref{fig2b}). If the movement rate of the infected people is much smaller (i.e. $\sigma\to 0$), the infected people may concentrate on the points of the highest risk.

Finally, when the infection rate is not directly proportional to the number of infected individuals (i.e., $0<p<1$), Theorems \ref{TH2.1}{\rm (ii)}, \ref{TH2.3}{\rm (ii)}, and \ref{TH2.4}{\rm (ii)} imply that implementing restrictions on the movement of susceptible and/or infected individuals may not help to eradicate the disease at any location within the domain. Moreover, there is a high density of the infected population in the highest-risk sites. Further studies are needed to determine whether these measures may help reduce the density of infected individuals.

In our companion paper (\cite{PSW2024b}), we will conduct a detailed analysis of the spatial profiles of model \eqref{model-2} with respect to small $d_S$ and/or $d_I$. Our findings will reveal significant differences between the results obtained from model \eqref{model-1} and model \eqref{model-2}. The comparisons of these contrasting results will be discussed in \cite{PSW2024b}.


\vskip10pt

\noindent\textbf{Acknowledgment}. The authors would like to thank the referees for value suggestions that lead to improvement of the manuscript. 

\vskip25pt
\begin{center}
DATA AVAILABILITY STATEMENT
\end{center}
Data sharing is not applicable to this article as no new data were created or analysed in this study.

\FloatBarrier

\vskip25pt
\bibliographystyle{plain}
%


\end{document}